\pdfoutput=1
\documentclass[preprint,12pt]{elsarticle}

\usepackage{amssymb}
\usepackage{algpseudocode}
\usepackage{algorithm}
\usepackage{tablefootnote}
\usepackage{amsmath}
\usepackage{placeins}
\usepackage{booktabs}
\usepackage{siunitx}
\usepackage{natbib}
\usepackage{subcaption}
\usepackage{mathrsfs}
\usepackage{hyperref}

\usepackage[table]{xcolor} 
\usepackage{enumitem}

\usepackage{amsmath}
\usepackage{pgfplots}
\usepackage[left=2.5cm,right=2.5cm,top=2.5cm,bottom=2.5cm]{geometry}
\usetikzlibrary{patterns,calc}
\pgfplotsset{compat=newest}
\usepgfplotslibrary{groupplots}

\pgfplotscreateplotcyclelist{customcolorlist}{
    {red, solid, mark=o,thick, draw opacity=1, fill opacity=1},           
    {red, solid, mark=o,thick, draw opacity=0.5, fill opacity=0.5}, 
    {blue, solid, mark=square,thick, draw opacity=1, fill opacity=1},          
    {blue, solid, mark=square,thick, draw opacity=0.5, fill opacity=0.5}, 
    {green!60!black, solid, mark=diamond,thick, draw opacity=1, fill opacity=1}, 
    {green!60!black, solid, mark=diamond,thick, draw opacity=0.5, fill opacity=0.5},  
    {orange!60!black, solid, mark=triangle,thick, draw opacity=1, fill opacity=1}, 
    {orange!60!black, solid, mark=triangle,thick, draw opacity=0.5, fill opacity=0.5}  
}

\pgfplotsset{
    grid=both,
    grid style={dashed, gray!30},
    ylabel near ticks,
    xlabel near ticks,
    enlargelimits,
	legend style={fill=white, fill opacity=0.6, draw opacity=1,draw=black,text opacity=1, font=\footnotesize},
    legend cell align={left},
    legend image with text/.style={
    legend image code/.code={%
        \node[anchor=center] at (0.3cm,0cm) {#1};
    }
},
    legend pos = north west,
tick label style={font=\footnotesize},
label style={font=\footnotesize},
}

\begin{document}

\begin{frontmatter}



\title{IsoGeometric Suitable Coupling Methods for Partitioned Multiphysics Simulation with Application to Fluid-Structure Interaction} 


\author[labela]{Jing-Ya Li} 
\author[labelb]{Hugo M. Verhelst}
\author[labelc]{Henk den Besten}
\author[labela]{Matthias M\"{o}ller}
\affiliation[labela]{organization={Delft University of Technology, Delft Institute of Applied Mathematics},
            addressline={Mekelweg 4}, 
            city={Delft},
            postcode={2628CD}, 
            country={The Netherlands}}

\affiliation[labelb]{organization={University of Pavia, Department of Civil Engineering and Architecture},
            addressline={ Via Adolfo Ferrata 3}, 
            city={Pavia},
            postcode={27100 PV}, 
            country={Italy}}

\affiliation[labelc]{organization={Delft University of Technology, Maritime and Transport Technology Department},
            addressline={Mekelweg 2}, 
            city={Delft},
            postcode={2628CD}, 
            country={The Netherlands}}

\begin{abstract}This paper presents spline-based coupling methods for partitioned multiphysics simulations, specifically designed for isogeometric analysis (IGA) based solvers. Traditional vertex-based coupling approaches face significant challenges when applied to IGA solvers, including geometric accuracy issues, interpolation errors, and substantial communication overhead. The methodology draws on the IGA mathematical framework to deliver coupling solutions that preserve the high-order continuity and exact geometric representation of splines. We develop two complementary strategies: (1) a spline-vertex coupling method that enables efficient interaction between IGA and conventional solvers, and (2) a fully isogeometric coupling approach that maximizes accuracy for IGA-to-IGA communication. 

Both theoretical analysis and extensive numerical experiments demonstrate that our spline-based methods significantly reduce communication overhead compared to traditional approaches while simultaneously enhancing geometric accuracy through exact boundary representation and maintaining higher-order solution continuity across the coupled interfaces.
We quantitatively confirm the communication efficiency benefits through systematic measurements of both transfer times and data volumes across various mesh refinement levels, with experimental results closely aligning with our theoretical predictions.
Our benchmark studies further demonstrate the geometric fidelity advantages through exact boundary representation, while also highlighting how the inherent mathematical structure of splines naturally preserves solution derivatives across interfaces without requiring additional computation or specialized transfer algorithms. This work not only provides efficient coupling strategies tailored to IGA-based solvers but also establishes a practical bridge between IGA and traditional discretization methods in partitioned multiphysics simulations. By offering viable options for coupling conventional solvers with IGA-based components, our approach enables broader adoption of IGA in established simulation workflows while ensuring accurate and high-performance interface communications.
\end{abstract}




\end{frontmatter}



\section{Introduction}

Fluid–structure interaction (FSI) remains a fundamental challenge in multiphysics simulations, governing a wide range of engineering applications, such as the aerodynamic response of aircraft wings, the pulsatile behavior of cardiovascular systems, and the structural resilience of bridges and offshore structures under hydrodynamic loading. These problems exhibit inherently bidirectional coupling, where the fluid exerts pressure on the structure, causing deformations that, in turn, modify the fluid domain, creating a nonlinear and dynamic feedback loop. Ensuring numerical stability, accuracy, and efficiency in such simulations is especially challenging when heterogeneous discretizations are utilized for the fluid and structural solvers, requiring robust interface coupling approaches \cite{NHoster_spline_based_FEM,NHoster_NURBS_based_coupling, FSI_Book}.

\subsection{Partitioned Methods for Multiphysics Problems}

FSI problems are typically tackled via one of the two following solution strategies: monolithic and partitioned formulations. The monolithic approach integrates fluid and structural solvers into a single, coupled system, yielding superior numerical stability and convergence behavior \cite{Monolisthic_vs_partitioned} but at the cost of flexibility, often necessitating tightly integrated code bases and customized solvers. In contrast, partitioned methods - central to the present study - decouple the problem into distinct subdomains, allowing fluid and structural equations to be solved independently and with different solvers \cite{benelhadj2019partitioned}. These subproblems communicate through the iterative exchange of boundary data across the interface. The modularity provided by partitioned schemes promotes solver reuse, algorithmic specialization, and compatibility with legacy software \cite{The_FSI_paper,partitioned_fsi_yulong}. Moreover, partitioned formulations are often favored over monolithic ones because they offer the flexibility to couple pre-existing solvers, facilitating the integration of specialized codes developed for each physical domain \cite{Chourdakis2021}. However, this same modularity introduces significant challenges: ensuring numerical stability, enforcing kinematic and dynamic consistency, and preserving accuracy across potentially non-conforming meshes are nontrivial tasks \cite{X_Jiao_conservative,X_Jiao_refinement_based,Kim_interface_element}.

Partitioned methods involve two main coupling aspects: a temporal and spatial one. Temporal coupling controls how and when data are exchanged between the fluid and structural solvers. This can range from simple staggered approaches to more advanced methods, such as fixed-point or quasi-Newton iterations that improve stability and accuracy \cite{X_Jiao_temporal}. Spatial coupling, on the other hand, deals with how data is transferred across the interface — especially when the fluid and structural meshes do not match. This usually involves interpolation or projection between different meshes \cite{X_Jiao_conservative, partitioned_fsi_yulong, Rendall_mesh_rbf}. It is difficult to achieve both high accuracy and good computational performance at the same time. The conventional approach of increasing mesh resolution to improve accuracy directly conflicts with the need to minimize communication overhead, leading to diminishing returns as the size of the problem increases.

\subsection{Spline-based Coupling: A Paradigm Shift}

Isogeometric Analysis (IGA), pioneered by Hughes et al. \cite{HUGHES20054135}, presents a powerful alternative to traditional finite element approaches by unifying geometry representation and solution approximation. In IGA, the same mathematical basis functions — typically Non-Uniform Rational B-Splines (NURBS) or B-splines — are used to define both the computational domain and the solution field \cite{cottrell2009isogeometric}. This seamless integration enables exact geometric representation — even on coarse meshes — and supports higher-order continuity across elements: features especially advantageous for problems dominated by interface behavior. 

Partitioned FSI methods have gained widespread adoption over the past decades due to their flexibility in leveraging specialized solvers for different physical domains. However, a persistent challenge in this context lies in the accurate and efficient transfer of data across the fluid–structure interface, particularly when nonmatching spatial discretizations are involved \cite{MALENICA2014983}. Spline-based coupling techniques, especially those employing NURBS or B-splines, have shown significant potential to address these issues. For example, Hosters et al.\ \cite{NHoster_spline_based_FEM} first developed a NURBS-enhanced FSI framework using spline descriptions on both fluid and structure sides, enabling accurate and direct data exchange. Kamensky et al.\ \cite{Kamensky2015} introduced a stabilized immersogeometric formulation, showing improved stability properties when using collocated constraints across thin, deforming interfaces.

Further developments in spline technology have expanded the scope of isogeometric coupling. Floating IGA (FLIGA), introduced by Hille et al.\ \cite{Hille2022}, enables analysis under extreme deformation by adapting basis functions along deformation-dependent directions. Similarly, Rosa et al.\ \cite{Rosa2022} proposed a blended IGA-FEM approach for modeling fracture propagation without remeshing, further illustrating the adaptability of spline-based representations to evolving geometries. Guarino et al.\ \cite{Guarino2024} presented a novel penalty-based coupling scheme for trimmed non-conformal shell patches, relevant for industrial applications where surface continuity and physical correctness must be preserved.

As the complexity of FSI problems increases — particularly in scenarios involving large structural deformations or compressible flows — robust coupling strategies become even more critical. While spline-based methods offer smooth and accurate interface representations, recent advances in finite element-based coupling have tackled these challenges from a different perspective. Rajanna et al.\ \cite{Rajanna_NonMatching_FEM} proposed a compressible FSI framework using nonmatching finite element meshes, suitable for high-speed flow applications. Klöppel et al.\ \cite{Klöppel_NonMatching_FEM} developed a dual mortar formulation to improve force and displacement consistency across nonconforming FEM interfaces. These contributions are important within the FEM paradigm, though they differ fundamentally in how geometry and continuity are handled compared to spline-based methods.

Despite advances in interface accuracy and stability, the issue of computational efficiency remains unresolved — particularly in large-scale, strongly coupled simulations. Data interpolation between nonmatching meshes can become a performance bottleneck. Radial Basis Function (RBF) methods have gained traction as black-box interpolants for multiphysics applications, as seen in the work of Lindner et al.\ \cite{Lindner2017} and its parallelized extension by Schneider et al.\ \cite{Schneider2023}. However, while RBFs offer flexibility, they may lack the conservation and variational structure desired in high-fidelity IGA settings. De Boer et al.\ \cite{deBoer2008} highlighted that consistent interpolation often outperforms conservative approaches in terms of accuracy and stability — particularly for quasi-1D FSI problems. In contrast, spline-based coupling inherently supports both consistency and high-order smoothness through function-space projection, offering a more principled route for interface treatment.

In contrast to existing point-based and spline-based coupling methods that mainly interpolate interface data at discrete points, this paper proposes a fundamentally new approach for partitioned FSI simulations by directly expressing interface fields as smooth, continuous spline functions. The key innovation lies in directly utilizing the inherent advantages of IGA to represent the interface as a continuous spline geometry, rather than as isolated interpolation points. In the proposed framework, both the structural and fluid participants maintain their own spline-defined interface fields, allowing for a consistent and symmetric field-to-field coupling. This shift in paradigm offers several compelling advantages: it preserves geometric fidelity for complex curvilinear boundaries, reduces communication overhead via compact control point exchange, improves numerical stability due to high-order continuity, and inherently handles nonmatching meshes through direct value evaluation, eliminating the need for ad hoc interpolation schemes.

The proposed methodology treats the interface as a mathematically structured entity rather than a numerical artifact. This approach maintains physical consistency while achieving high accuracy with reduced computational cost, effectively addressing the longstanding accuracy-efficiency trade-off in spatial coupling for partitioned FSI systems.

\subsection{Research Contributions and Paper Structure}

This work introduces key contributions to multiphysics simulation. While the methods are broadly applicable to general multiphysics problems, we focus without loss of generality on fluid–structure interaction as a representative application throughout this study.
\begin{itemize}
    \item A theoretical framework for interface coupling in partitioned FSI based on IGA, specifically addressing challenges with non-matching meshes while maintaining accuracy and computational efficiency.
    \item Two coupling strategies: (i) a hybrid approach that combines IGA-based structural solvers with traditional fluid solvers, and (ii) a fully isogeometric method for high-accuracy IGA-to-IGA coupling.
    \item An analysis of communication overhead, showing that our spline-based method significantly reduces communication cost, especially with finer meshes.
    \item Verification and validation using standard FSI benchmarks, demonstrating a strong balance between accuracy and efficiency.
    \item Implementation as an open-source software package integrated into the preCICE multiphysics framework, featuring a modular plug-and-play design compatible.
\end{itemize}

The remainder of the paper is organized as follows. Section~\ref{Section:Governing} introduces the governing equations used in this study. Section~\ref{Section:IGA} provides an overview of IGA, which serves as the spatial discretization method for the structural domain in this work. 
Section~\ref{Section:Numerical Methods} explains the numerical methods for both fluid and structural solvers and presents our spline-based coupling approach. Section~\ref{section:verification} provides verification and validation results to demonstrate the accuracy of our method. 
The concluding Section~\ref{section:conclusion} ends this study, summing up the main results and opening up for further research. 

The methodology presented in this work has been implemented in \texttt{gsPreCICE}, developed as part of this research. As the first dedicated IGA solver adapter for the preCICE (Precise Code Interaction Coupling Environment) \cite{Chourdakis2021} framework, \texttt{gsPreCICE} is open-source and officially integrated as a module within the \texttt{G+Smo} (Geometry + Simulation Modules)  library \cite{mantzaflaris2019overview,juttler2014geometry,gismo_paper}, which is a C++ library dedicated to IGA.





\section{Governing equations}
\label{Section:Governing}
Partitioned algorithms tackle the distinct challenges posed by fluid dynamics and structural mechanics by separately solving each domain and exchanging data at a shared interface. While such approaches are broadly applicable to a wide range of multiphysics problems, the scope here is limited to FSI, involving the coupling of incompressible fluid flow with elastic structural deformation in multidimensional spaces.

In this section, we briefly review the governing equations of incompressible fluid flow and elastic structural deformation used throughout the paper. Clearly, restating these equations here serves two purposes: (a) it establishes a consistent notation that is essential for accurately describing our numerical methods and coupling strategies, and (b) it provides context for subsequent discussions on interface conditions and numerical discretizations.

\subsection{Incompressible Flow}
The fluid dynamics within the domain are governed by the Navier-Stokes equations for incompressible flow, outlined as follows:

\begin{subequations}
\begin{align}
\rho^f \left( \frac{\partial\mathbf{u}^f}{\partial t} + (\mathbf{u}^f \cdot \nabla) \mathbf{u}^f \right) 
    &= - \nabla p^f + \nabla \cdot(\mu^f \nabla \mathbf{u}^f) + \mathbf{f}^f \quad \text{in } \Omega^f\times[0,T], \label{eq:momentum} \\
\nabla \cdot \mathbf{u}^f 
    &= 0 \quad \text{in } \Omega^f\times[0,T], \label{eq:continuity} \\
\mathbf{u}^f 
    &= \mathbf{g}^f \quad \text{on } \Gamma_{D}^f \times[0, T], \label{eq:Dirichlet} \\
(\mu^f \nabla \mathbf{u}^f - p^f \mathbf{I}) \cdot \mathbf{n}^f 
    &= \mathbf{h}^f \quad \text{on } \Gamma_N^f\times[0,T]. \label{eq:Neumann}
\end{align}
\end{subequations}
Within these equations, $\rho^f$, $\mu^f$, and $\mathbf{f}^f$ represent the fluid's density, dynamic viscosity, and external body force vector, respectively. $\Gamma_{D}^f$ and $\Gamma_{N}^f$ denote the Dirichlet and Neumann boundaries of the fluid domain $\Omega^f$. This set of equations, (\ref{eq:momentum}) through (\ref{eq:Neumann}), is utilized to solve for the fluid's velocity field $\mathbf{u}^f(\mathbf{x},t)$ and pressure $p^f(\mathbf{x},t)$ at any given time $t \in [0,T]$.

\subsection{Structure deformation}
The solid displacement field, $\mathbf{d}^s(\mathbf{x},t)$, which responds to external forces, is governed by Newton's second law of motion, as expressed below:

\begin{subequations} \label{eq:structure_full}
\begin{equation}
\frac{d^2 \mathbf{d}^s}{d t^2}-\frac{1}{\rho_0^s} \boldsymbol{\nabla}_0 \cdot\left(\mathbf{S}^s \mathbf{F}^T\right) = \mathbf{b}^s \quad \text{in } \Omega_0^s \times [0,T] \label{eq:structure_motion}
\end{equation}
\end{subequations}

Here, $\rho_0^s$ represents the solid's density and $\mathbf{b}^s$ denotes the external body force acting on it. This equation is derived from a total Lagrangian framework and is hence formulated with respect to the unstrained reference configuration $\Omega_0^s$. All quantities and operators in this formulation are referenced with a subscript '0' to indicate this baseline state. Unlike the Cauchy stress tensor $\mathbf{T}^s$, the internal stress here is expressed using the second Piola-Kirchhoff stress tensor $\mathbf{S}^s = \operatorname{det}(\mathbf{F}) \mathbf{F}^{-1} \mathbf{T}^s \mathbf{F}^{-T}$, where $\mathbf{F}$ is the deformation gradient. A Hookean or the St. Venant-Kirchhoff material model typically forms the basis of the constitutive equation, contributing to the geometric nonlinearity of the solid mechanics problem. The problem formulation is completed with initial zero displacement conditions and appropriate boundary constraints on $\partial \Omega^s$. The condition reads:
\begin{subequations} \label{eq:structure_bc}
\begin{align}
\mathbf{d}^s &= \mathbf{g}^s \quad \text{on } \Gamma_{D}^s \times [0, T] \label{eq:structure_dirichlet} \\
\mathbf{F}\mathbf{S}^s\mathbf{n}_0 &= \mathbf{h}^s \quad \text{on } \Gamma_{N}^s \times [0, T] \label{eq:structure_neumann}
\end{align}
\end{subequations}

\subsection{Coupling conditions}
In the partitioned approach, the fluid problem and the structure problem are solved independently, potentially by different solvers. To ensure the conservation of mass, momentum, and mechanical energy over the shared interface $\Gamma^{fs} = \Gamma^f\cap\Gamma^s$, the solution along the interface for the two subproblems has to satisfy kinematic continuity:
\begin{subequations}
    \begin{equation}
        \mathbf{d}^f = \mathbf{d}^s \quad \text{on } \Gamma^{fs}\times[0,T]
    \end{equation}
    \begin{equation}
        \mathbf{u}^f = \mathbf{u}^s \quad \text{on } \Gamma^{fs}\times[0,T]
    \end{equation}
\end{subequations}
where $\mathbf{d}^f$ and $\mathbf{u}^f$ represent the displacement and velocity of the moving fluid boundary, while $\mathbf{d}^s$ and $\mathbf{u}^s$ are those of the structure.

The implied accelerations are equivalent too. And the dynamic continuity for the equilibrium of stress:
\begin{equation}
    \mathbf{T}^f\cdot\mathbf{n}^f = -\mathbf{T}^s\cdot\mathbf{n}^s \quad \text{on } \Gamma^{fs}\times[0,T]
\end{equation}
where $\mathbf{n}^f$ and $\mathbf{n}^s$ are the interface unit normal vectors pointing outwards from the corresponding domains. Note that $\mathbf{n}^f = -\mathbf{n}^s$ on the interface. These conditions ensure mass, momentum, and energy conservation across the interface \cite{FSI_Book}.

\section{Isogeometric Analysis Overview}
\label{Section:IGA}
In this work, we adopt a partitioned coupling approach, treating fluid and structure as separate subproblems that are solved independently. For the structural subproblem, we employ IGA, which offers both geometric exactness and higher continuity in the solution space.

IGA, introduced by Hughes et al. \cite{HUGHES20054135}, enhances geometric accuracy and promotes closer integration between CAD models and numerical analysis by using the same basis functions for geometry representation and solution approximation. Typically, CAD geometries are represented by NURBS or b-splines, defined by a set of control points and associated rational basis functions. A NURBS surface can be mathematically expressed as:
\begin{equation}
\overline{\mathbf{x}}(\boldsymbol{\xi})=\sum_{i=1}^{n_c} R_i^p(\boldsymbol{\xi}) \mathbf{x}_i,\quad \text{with }\boldsymbol{\xi}=(\xi^1,\xi^2),
\end{equation}
where $n_c$ is the number of control points, $\mathbf{x}_i \in \mathbb{R}^{d}$ represents the control points, $R_i^p(\boldsymbol{\xi})$ are rational basis functions of polynomial degree $p$ defined on the parameter domain $\boldsymbol{\Xi}\subset\mathbb{R}^d$, and $d\in\{2,3\}$ is the spatial dimension.

IGA provides an ideal framework for our coupling methods because it offers smooth and continuous representations of solution fields across the interface, which are critical for stable and accurate multiphysics coupling. In addition, it enables precise geometric descriptions, which help avoid approximation errors in curved or evolving interfaces.

\section{Numerical Methods}
\label{Section:Numerical Methods}
In partitioned FSI simulations, the fluid and structure subproblems are solved using separate solvers, each tailored to its respective physical field. This modular strategy offers considerable flexibility: established solvers can be reused or further developed independently, without requiring changes to the overall framework \cite{Chourdakis2021}. 

At the interface between fluid and structure, the two solvers often use different meshes — and sometimes even different types of numerical methods. Because of this, we need a reliable way to transfer data (like forces and displacements) between them in a way that is accurate and physically consistent. One of the main contributions of this paper is a spline-based coupling method, which is especially effective when combined with IGA. This method helps maintain smooth and consistent data transfer across nonmatching meshes.

Building on this motivation, the remainder of this section focuses on spatial coupling strategies that bridge the fluid and structural solvers, with particular emphasis on the proposed spline-based methods that support the consistency and efficiency of the partitioned framework.

\FloatBarrier
\subsection{Evolution of Spatial Coupling Schemes: From Vertex-based to Spline-based Methods}
Spatial coupling schemes are critical components in partitioned FSI simulations, determining the manner and accuracy of data exchange at the interface. We consider three representative approaches: vertex-vertex, spline-vertex, and spline-spline coupling, that range from conventional point-based methods to the proposed spline-based formulations.

\subsubsection{Vertex-based Coupling Method: Traditional Approach and Its Limitations}
The vertex-based coupling method facilitates communication between an IGA-based structural solver and traditional solvers based on discretization techniques such as finite volume (FV) or finite element methods (FEM). In this approach, the data exchanged between solvers is defined at discrete points along the interface — typically mesh vertices or quadrature points.
\FloatBarrier
\begin{figure}[h!]
    \centering
    \includegraphics[width=0.9\linewidth]{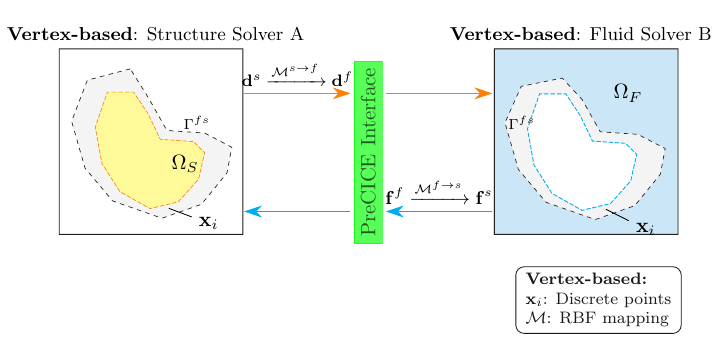}
    \caption{Illustration of vertex-based solvers coupled with IGA-based solvers, with information exchanged at quadrature points/vertices.}
    \label{fig:quad-coupling}
\end{figure}
\FloatBarrier
The vertex-based coupling procedure is performed in three main steps:
\begin{enumerate}[label=(\alph*)]
\item \textbf{Data mapping for the discrete forces:} Vertex-based coupling requires discrete stress components at each fluid node along the FSI interface $\Gamma^{fs}$. To transfer these stresses between nonmatching meshes, a range of mesh-to-mesh interpolation techniques are available, including nearest-neighbor, linear projection, and mortar methods. In this work, we use RBF interpolation, implemented via the preCICE library \cite{Schneider2023, Lindner2017}, due to its flexibility, smoothness, and robustness for scattered, non-conforming interface meshes. Compared to nearest-neighbor or linear methods, RBF interpolation offers higher accuracy and smoother field transfer without requiring mesh alignment or auxiliary mesh constructions \cite{Schneider2023}. This enables accurate interpolation of stress values from the fluid mesh $\Gamma^f$ to the structural mesh $\Gamma^s$.

Each interface mesh is defined by a set of vertices $\{ \mathbf{x}_i \in \mathbb{R}^d \, | \, i = 1, \dots, N \}$ in $d \in \{2,3\}$ spatial dimensions. Let $\{\boldsymbol{\sigma}_i \in \mathbb{R}^m \, | \, i = 1, \dots, N\}$ denote the stress tensor or vector values (with $m = d$ or $m = d(d+1)/2$ depending on tensor representation) at the fluid mesh vertices. For RBF interpolation, each stress component is treated as a scalar field. Given the $k$-th stress component $\{\sigma_{i,k} \in \mathbb{R} \, | \, i = 1, \dots, N\}$ at the input vertices, the interpolant $\mathcal{I}_{\sigma_k} : \mathbb{R}^d \to \mathbb{R}$ for the $k$-th component is defined as:
\begin{equation}
\mathcal{I}_{\sigma_k}(\mathbf{x}) = \sum_{i=1}^N \lambda_{i,k} \, \phi\left(\left\| \mathbf{x} - \mathbf{x}_i \right\| \right) + \beta_{0,k} + \boldsymbol{\beta}_{l,k}^T \mathbf{x},
\end{equation}
where $N$ is the number of input vertices, $\lambda_{i,k} \in \mathbb{R}$ are the RBF coefficients, $\phi : \mathbb{R} \to \mathbb{R}$ is the radial basis function, and $\beta_{0,k} \in \mathbb{R}$, $\boldsymbol{\beta}_{l,k} \in \mathbb{R}^d$ are coefficients of a linear polynomial ensuring consistency for linear stress fields. The interpolation conditions are enforced as:
\begin{equation}
\mathcal{I}_{\sigma_k}(\mathbf{x}_i) = \sigma_{i,k}, \quad i = 1, \dots, N,
\end{equation}
with regularization constraints to guarantee uniqueness:
\begin{equation}
\sum_{i=1}^N \lambda_{i,k} = 0 \quad \text{and} \quad \sum_{i=1}^N \lambda_{i,k} \mathbf{x}_i = \mathbf{0}.
\end{equation}

\item \textbf{Transfer of forces onto the structure:} The process is repeated for all stress components via preCICE. The interpolant $\mathcal{I}$ can yield a matrix $\mathcal{M}^{f\to s}$ defining a linear mapping of the input forces $\mathbf{f}^f$ from the fluid mesh to forces $\mathbf{f}^s$ in the structural mesh:
\begin{equation}
    \mathbf{f}^s = \mathcal{M}^{f\to s}\mathbf{f}^f 
    \quad \text{with} \quad 
    \mathbf{f}^f= \begin{pmatrix}
        \mathbf{f}^f(\mathbf{x}_1) \\
        \mathbf{f}^f(\mathbf{x}_2) \\
        \vdots \\
        \mathbf{f}^f(\mathbf{x}_N)
    \end{pmatrix}
\label{eq:force-interpolation-vertex}
\end{equation}
When the row-sum of $\mathcal{M}^{f\to s}$ is equal to one, the interpolation is consistent, which means that constant values are interpolated exactly. On the other hand, when the column-sum of $\mathcal{M}^{f\to s}$ equals one the interpolation is conservative.

\item \textbf{Transfer of deformation onto the fluid:} The computed structural deformation has to be transferred back to the fluid in order to account for the mesh deformation. Once again, this procedure is straightforward due to the use of the RBF method in preCICE. The deformation is computed at every quadrature point of the structural element:
\begin{equation}
    \mathbf{d}^f = \mathcal{M}^{s\to f}\mathbf{d}^s \quad \text{with} \quad \mathbf{d}^s= \begin{pmatrix}
\mathbf{d}^s(\mathbf{x}_1) \\
\mathbf{d}^s(\mathbf{x}_2) \\
\vdots \\
\mathbf{d}^s(\mathbf{x}_n)
\end{pmatrix}
\end{equation}
where $\mathbf{d}^s$ represents the structural displacements at discrete structural points $\mathbf{x}_i$, and $\mathbf{d}^f$ represents the fluid mesh displacements at discrete fluid mesh points. $\mathcal{M}^{s\to f}$ represents the interpolation operator, which maps the structural displacement onto the fluid mesh on the IGA mesh cell level.
\end{enumerate}

Although vertex-based methods are widely used due to their simplicity and compatibility with standard finite element frameworks, they face several significant challenges when dealing with complex interface dynamics \cite{X_Jiao_refinement_based,Kim_interface_element}:
\begin{enumerate}
    \item Geometric accuracy issues, especially for complex or curved interfaces \cite{N_Hoster_quasi_newton}.
    \item Interpolation errors that can lead to numerical inaccuracies in force and displacement transfer processes.
    \item Significant increase in communication overhead as interface meshes are refined.
\end{enumerate}

These limitations motivate the development of more accurate and efficient spline-based coupling methods, which are introduced in this paper.

\subsection{Spline-based Coupling Method}
\label{subsec:spline-based-coupling}

The spline-based coupling method leverages analytical spline functions to represent the interface, providing a more efficient and accurate approach for data exchange between solvers. Unlike vertex-based methods that operate on discrete points, spline-based coupling utilizes a continuous mathematical description of both interface geometry and field variables.

We propose two types of spline-based coupling methods:  
\begin{itemize}
    \item Spline-vertex coupling method: one solver is IGA-based.
    \item IGA-IGA coupling method: both solvers are IGA-based.
\end{itemize}

\subsubsection{Spline-vertex Coupling Method}
\label{subsec:spline-vertex-based-coupling}
We first consider the spline-vertex coupling method, which is designed for scenarios in which only one side of the interface—typically the structure—is modeled using IGA, while the opposing side (e.g., the fluid) employs a standard finite element or finite volume discretization. In this case, the interface geometry is represented using a spline on the IGA side, and data is transferred between the spline surface and the vertex points of the non-IGA mesh using projection or interpolation techniques. 
\FloatBarrier
\begin{figure}[htbp!]
    \centering
    \includegraphics[width=0.9\linewidth]{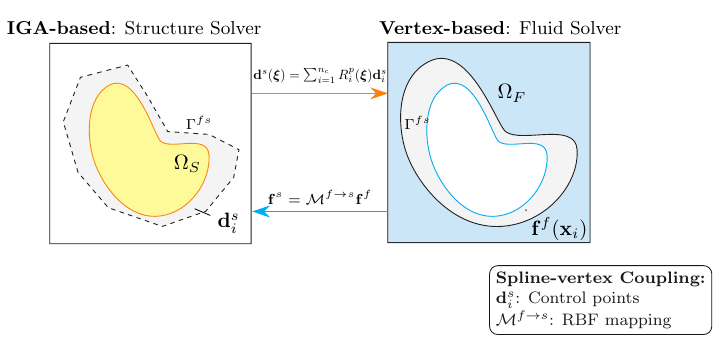}
    \caption{Spline-based coupling between IGA-based structural solver and vertex-based fluid solver, with interface data exchanged via continuous spline representation.}
    \label{fig:spline-coupling}
\end{figure}
\FloatBarrier
The coupling procedure consists of two key steps:
\begin{itemize}
\item \textbf{(a) Force transfer (fluid to structure):} Force data from the fluid solver is transferred to the structure using a RBF interpolation method:
\begin{equation}
    \mathbf{f}^s = \mathcal{M}^{f\to s}\mathbf{f}^f 
    \quad \text{with} \quad 
    \mathbf{f}^f= \begin{pmatrix}
        \mathbf{f}^f(\mathbf{x}_1) \\
        \mathbf{f}^f(\mathbf{x}_2) \\
        \vdots \\
        \mathbf{f}^f(\mathbf{x}_N)
    \end{pmatrix}
\label{eq:force-interpolation-vertex}
\end{equation}
This step maps the discrete fluid force data to the structural solver's representation.

\item \textbf{(b) Displacement transfer (structure to fluid):} The structural displacement field, computed in terms of control point displacements $\mathbf{d}_i^s$, is evaluated at any interface point using:
\begin{equation}
\mathbf{d}^s(\boldsymbol{\xi}) = \sum_{i=1}^{n_c} R_i^p(\boldsymbol{\xi}) \mathbf{d}_i^s
\end{equation}
The fluid solver samples this continuous field at required interface points to update its mesh.
\end{itemize}


\subsubsection{IGA-IGA Communication: Optimal Spline-based Coupling}
\label{subsec:iga-iga-communication}

When both fluid and structural domains use isogeometric analysis or support spline representation (of solution/geometry), the IGA-IGA communication approach fully leverages spline representations on both sides of the interface, representing the most efficient coupling strategy within the spline-based family.

\begin{figure}[h!]
    \centering
    \includegraphics[width=0.9\linewidth]{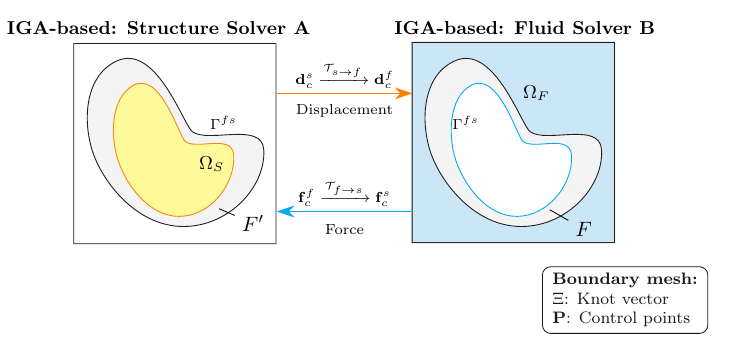}
    \caption{IGA-IGA communication with spline representations on both fluid and structural domains, maximizing the benefits of isogeometric analysis.}
    \label{fig:iga-iga-coupling}
\end{figure}

The communication procedure involves:
\begin{itemize}
\item \textbf{(a) Interface parameterization:} Ideally, both domains share a common parametric representation, enabling direct exchange of control point data. When this is not possible, a reparameterization procedure establishes correspondence between the two spline spaces by exchanging the knot vector along the interface. 
\item \textbf{(b) Stress transfer via control points:} Forces are directly expressed in terms of control point values. With matching parameterizations, they transfer directly; otherwise, a projection between spline spaces is performed:
\begin{equation}
\mathbf{f}^s_c = \mathcal{T}_{f \rightarrow s} \mathbf{f}^f_c
\end{equation}
where $\mathbf{f}^s_c$ and $\mathbf{f}^f_c$ are the control point force vectors for the structural and fluid domains, respectively, and $\mathcal{T}_{f \rightarrow s}$ is the transformation operator mapping between the two spline spaces.
\item \textbf{(c) Displacement transfer via control points:} Similarly, structural displacements computed at control points transfer directly or through transformation:
\begin{equation}
\mathbf{d}^f_c = \mathcal{T}_{s \rightarrow f} \mathbf{d}^s_c
\end{equation}

where $\mathbf{d}^s_c$ and $\mathbf{d}^f_c$ are the control point displacement vectors for the structural and fluid domains, respectively, and $\mathcal{T}_{s \rightarrow f}$ is the transformation operator mapping between the two spline spaces.
\end{itemize}

This approach provides minimal communication overhead, exact geometry representation, preservation of higher-order solution properties, enhanced numerical stability, and efficient handling of complex geometries. We will further discuss the advantage of reduced communication overhead in Section \ref{section:verification}.


\section{Coupling Methods Verification and Validation}
\label{section:verification}

In this section, we present a comprehensive verification and validation of the spatial coupling methods described previously. The verification procedure employs both the constant load and nonlinear load test cases to assess the accuracy, efficiency, and stability of the proposed coupling approaches.

\begin{figure}[h!]
    \centering
    \begin{minipage}{0.4\textwidth}
        \centering
        \includegraphics[width=\textwidth]{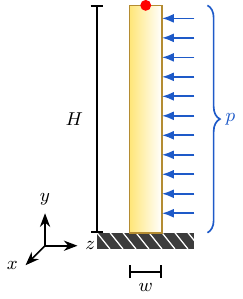}  
    \end{minipage}%
    \hfill
    \begin{minipage}{0.45\textwidth}
        \centering
        \renewcommand{\arraystretch}{1.2}  
        \begin{tabular}{@{}l l@{}}
            \toprule
            \textbf{Parameter} & \textbf{Value} \\
            \midrule
            Width ($W$)       & 0.50 m \\
            Height ($H$)       & 1.00 m \\
            Thickness ($t$)        & 0.05 m \\
            Density ($\rho$) & \si{3000 \text{ }kg/m^3} \\
            Distributed Load ($p$) & $5\times10^3 \mathbf{e}_z$ \si{Pa}  \\
            \bottomrule
        \end{tabular}
    \end{minipage}
    \caption{Constant load test setup for the vertical beam.}
    \label{figure:constant-vertical-beam-setup}
\end{figure}

We first consider a vertical beam that is clamped at its base and subjected to a constant distributed load along its free edge. The setup ensures that the force transmission is uniform and can be analytically validated.
\FloatBarrier
The problem setup, including geometry and loading conditions, is shown in Figure~\ref{figure:constant-vertical-beam-setup}. We compare the tip displacement computed by the standalone structural solver (reference solution) to the coupled solver solutions using vertex-based and spline-based coupling methods. 
\FloatBarrier

The structural solver computes the initial beam deformation under the applied load. The force distribution is mapped to the fluid mesh using the vertex-based and spline-based methods described earlier. The reconstructed forces are transferred back to the structure to verify accuracy in communication. Finally, the resulting tip displacement is compared with the theoretical solution, which is the monolithic application of the force along the beam within one solver.

\FloatBarrier
\begin{figure}[h!]
    \centering
    \begin{minipage}[t]{0.49\textwidth}
        \vspace{0pt}
        \begin{tikzpicture}
        \begin{groupplot}[
            group style={
                group name=myplot,
                group size=1 by 2,
                vertical sep=15pt,
                xticklabels at=edge bottom
            },
            width=\linewidth,
            height=4.5cm,
            grid=major,
            xlabel style={font=\small},
            ylabel style={font=\small},
            ticklabel style={font=\footnotesize},
            legend style={
                font=\footnotesize,
                draw=none,
                fill=none,
                at={(0.5,-1.55)},
                anchor=north,
                legend columns=2,
                column sep=10pt
            },
            legend cell align={left}
        ]

        \nextgroupplot[ylabel={$u_y\;[\mathrm{m}]$}]
        \addplot+[mark=square*, mark options={fill=blue}, smooth, mark repeat=10, line width=0.8pt] 
            table[x index=3, y index=1, col sep=comma]{validation_data/ConstantLoad/pointData-vertex-vertex.csv};
        \addlegendentry{Vertex-Vertex}

        \addplot+[mark=o, mark options={fill=red}, smooth, mark repeat=10, line width=0.8pt] 
            table[x index=3, y index=1, col sep=comma]{validation_data/ConstantLoad/pointData-IGA-vertex.csv};
        \addlegendentry{Spline-Vertex}

        \addplot+[mark=triangle*, mark options={fill=brown}, smooth, mark repeat=10, line width=0.8pt] 
            table[x index=3, y index=1, col sep=comma]{validation_data/ConstantLoad/pointData-IGA-IGA.csv};
        \addlegendentry{Spline-Spline}

        \addplot+[
            color=green,
            dashed,
            mark=*,
            mark options={fill=green},
            smooth,
            mark repeat=10,
            line width=0.8pt
        ] 
        table[x index=3, y index=1, col sep=comma]{validation_data/ConstantLoad/pointData-reference.csv};
        \addlegendentry{Reference}

        \nextgroupplot[xlabel={$t\;[\mathrm{s}]$}, ylabel={$u_z\;[\mathrm{m}]$}]
        \addplot+[mark=square*, mark options={fill=blue}, smooth, mark repeat=10, line width=0.8pt] 
            table[x index=3, y index=2, col sep=comma]{validation_data/ConstantLoad/pointData-vertex-vertex.csv};

        \addplot+[mark=o, mark options={fill=red}, smooth, mark repeat=10, line width=0.8pt] 
            table[x index=3, y index=2, col sep=comma]{validation_data/ConstantLoad/pointData-IGA-vertex.csv};

        \addplot+[mark=triangle*, mark options={fill=brown}, smooth, mark repeat=10, line width=0.8pt] 
            table[x index=3, y index=2, col sep=comma]{validation_data/ConstantLoad/pointData-IGA-IGA.csv};

        \addplot+[
            color=green,
            dashed,
            mark=*,
            mark options={fill=green},
            smooth,
            mark repeat=10,
            line width=0.8pt
        ] 
        table[x index=3, y index=2, col sep=comma]{validation_data/ConstantLoad/pointData-reference.csv};

        \end{groupplot}
        \end{tikzpicture}
    \end{minipage}%
    \hfill
    \begin{minipage}[t]{0.49\textwidth}
        \vspace{0pt}
        \begin{tikzpicture}
        \begin{groupplot}[
            group style={
                group size=1 by 2,
                vertical sep=15pt,
                xticklabels at=edge bottom
            },
            width=\linewidth,
            height=4.37cm,
            grid=major,
            xlabel style={font=\small},
            ylabel style={font=\small},
            ticklabel style={font=\footnotesize},
            legend style={
                at={(0.5,-1.55)},
                anchor=north,
                font=\footnotesize,
                draw=none,
                fill=none,
                legend columns=2,
                column sep=10pt
            },
            legend cell align={left}
        ]

        \nextgroupplot[
            ylabel={$\|u_y - \hat{u}_y\|\;[\mathrm{m}]$}
        ]
        \addplot+[mark=square*, mark options={fill=blue}, smooth, mark repeat=10, line width=0.8pt]
            table[x index=3, y index=4, col sep=comma]{validation_data/ConstantLoad/pointData-vertex-vertex.csv};
        \addlegendentry{Vertex-Vertex}

        \addplot+[mark=o, mark options={fill=red}, smooth, mark repeat=10, line width=0.8pt]
            table[x index=3, y index=4, col sep=comma]{validation_data/ConstantLoad/pointData-IGA-vertex.csv};
        \addlegendentry{Spline-Vertex}

        \addplot+[mark=triangle*, mark options={fill=brown}, smooth, mark repeat=10, line width=0.8pt]
            table[x index=3, y index=4, col sep=comma]{validation_data/ConstantLoad/pointData-IGA-IGA.csv};
        \addlegendentry{Spline-Spline}

        \nextgroupplot[
            xlabel={$t\;[\mathrm{s}]$},
            ylabel={$\|u_z - \hat{u}_z\|\;[\mathrm{m}]$}
        ]
        \addplot+[mark=square*, mark options={fill=blue}, smooth, mark repeat=10, line width=0.8pt]
            table[x index=3, y index=5, col sep=comma]{validation_data/ConstantLoad/pointData-vertex-vertex.csv};

        \addplot+[mark=o, mark options={fill=red}, smooth, mark repeat=10, line width=0.8pt]
            table[x index=3, y index=5, col sep=comma]{validation_data/ConstantLoad/pointData-IGA-vertex.csv};

        \addplot+[mark=triangle*, mark options={fill=brown}, smooth, mark repeat=10, line width=0.8pt]
            table[x index=3, y index=5, col sep=comma]{validation_data/ConstantLoad/pointData-IGA-IGA.csv};

        \end{groupplot}

        \node[font=\footnotesize] at (current bounding box.south) 
            [yshift=-0.8cm] {All results remain below machine precision ($\sim 10^{-15}$)};
        \end{tikzpicture}
    \end{minipage}
    \caption{Left: Comparison of tip displacement between reference solution and different coupling methods. Right: Relative $L^2$ error of different coupling methods compared to the reference solution.}
    \label{figure:constant-vertical-beam-tip-displacement}
\end{figure}
\FloatBarrier
Figure~\ref{figure:constant-vertical-beam-tip-displacement} compares the vertical beam tip displacements obtained via different coupling strategies against the monolithic reference solution (left), and the corresponding $L^2$-norm errors in the $y$- and $z$-displacement directions (right). All spline-based methods (spline-vertex and spline-spline) yield errors below machine precision ($\sim 10^{-15}$), indicating perfect agreement with the reference. The vertex-vertex method achieves the same level of accuracy in this constant load scenario. This demonstrates the consistency of all tested coupling approaches.

To further validate the accuracy and robustness of the coupling methods, we designed a nonlinear load test case. In this scenario, the vertical beam is subjected to a nonlinear load $\hat{p}(\xi) = -\xi\mathbf{e}_z$ that varies with parameter $\xi$. 

\FloatBarrier
\begin{figure}[h!]
    \centering
    \begin{minipage}{0.40\textwidth}
        \centering
        \includegraphics[width=\textwidth]{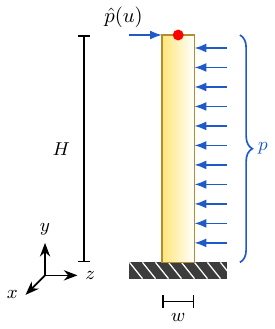}  
    \end{minipage}%
    \hfill
    \begin{minipage}{0.45\textwidth}
        \centering
        \renewcommand{\arraystretch}{1.2}  
        \begin{tabular}{@{}l l@{}}
            \toprule
            \textbf{Parameter} & \textbf{Value} \\
            \midrule
            Width ($W$)       & 0.50 m \\
            Height ($H$)       & 1.00 m \\
            Thickness ($t$)        & 0.05 m \\
            Density ($\rho$) & \si{3000\text{ } kg/m^3} \\
            Distributed Load ($p$) & $5\times10^3 \mathbf{e}_z$ \si{Pa}  \\
            Nonlinear Load ($\hat{p}(\xi)$) & $-\xi\mathbf{e}_z$\\
            \bottomrule
        \end{tabular}
    \end{minipage}
    \caption{Nonlinear load test setup for the vertical beam. The nonlinear load $\hat{p}(\xi)$ varies linearly with parameter $\xi$, where $\xi$ is the normalized coordinate along the beam height, with $\xi=0$ at the fixed bottom end and $\xi=1$ at the free top end.}
    \label{figure:nonlinear-vertical-beam-setup}
\end{figure}

Figure~\ref{figure:nonlinear-vertical-beam-setup} illustrates the nonlinear load test setup. In this test, we compare the reference solution (direct calculation by a monolithic solver) with solutions obtained using different coupling methods (vertex-based and spline-based) under identical load conditions. 
\FloatBarrier

\FloatBarrier
\begin{figure}[h!]
\centering
\begin{tikzpicture}[]
\begin{groupplot}[
    group style={
        group name=myplot,
        group size=1 by 2,
        vertical sep=20pt,
        xticklabels at=edge bottom
    },
    width=13cm,
    height=5.5cm,
    grid=major,
    xlabel style={font=\small},
    ylabel style={font=\small},
    ticklabel style={font=\footnotesize},
    legend style={
        at={(0.5,-1.5)},
        anchor=north,
        legend columns=3,
        font=\footnotesize,
        draw=none,
        fill=none,
        /tikz/every even column/.append style={column sep=0.5cm}
    },
    legend cell align=center
]

\nextgroupplot[ylabel={Displacement $u_y\;[\mathrm{m}]$}]
\addplot+[mark=square*, mark options={fill=gray}, smooth, mark repeat=5, line width=1pt] 
    table[x index=3, y index=1, col sep=comma] {validation_data/NonlinearLoad/pointData-vertex-vertex.csv};
\addlegendentry{Vertex-Vertex}

\addplot+[mark=o, mark options={fill=blue}, smooth, mark repeat=5, line width=1pt] 
    table[x index=3, y index=1, col sep=comma] {validation_data/NonlinearLoad/pointData-IGA-vertex.csv};
\addlegendentry{Spline-Vertex}

\addplot+[mark=triangle*, mark options={fill=red}, smooth, mark repeat=5, line width=1pt] 
    table[x index=3, y index=1, col sep=comma] {validation_data/NonlinearLoad/pointData-IGA-IGA.csv};
\addlegendentry{Spline-Spline}

\nextgroupplot[xlabel={$t\;[\mathrm{s}]$}, ylabel={Displacement $u_z\;[\mathrm{m}]$}]
\addplot+[mark=square*, mark options={fill=gray}, smooth, mark repeat=5, line width=1pt] 
    table[x index=3, y index=2, col sep=comma] {validation_data/NonlinearLoad/pointData-vertex-vertex.csv};

\addplot+[mark=o, mark options={fill=blue}, smooth, mark repeat=5, line width=1pt] 
    table[x index=3, y index=2, col sep=comma] {validation_data/NonlinearLoad/pointData-IGA-vertex.csv};

\addplot+[mark=triangle*, mark options={fill=red}, smooth, mark repeat=5, line width=1pt] 
    table[x index=3, y index=2, col sep=comma] {validation_data/NonlinearLoad/pointData-IGA-IGA.csv};
\end{groupplot}
\end{tikzpicture}
\caption{Tip displacement versus time under nonlinear loading.}
\label{figure:nonlinear-vertical-beam-tip-displacement}
\end{figure}
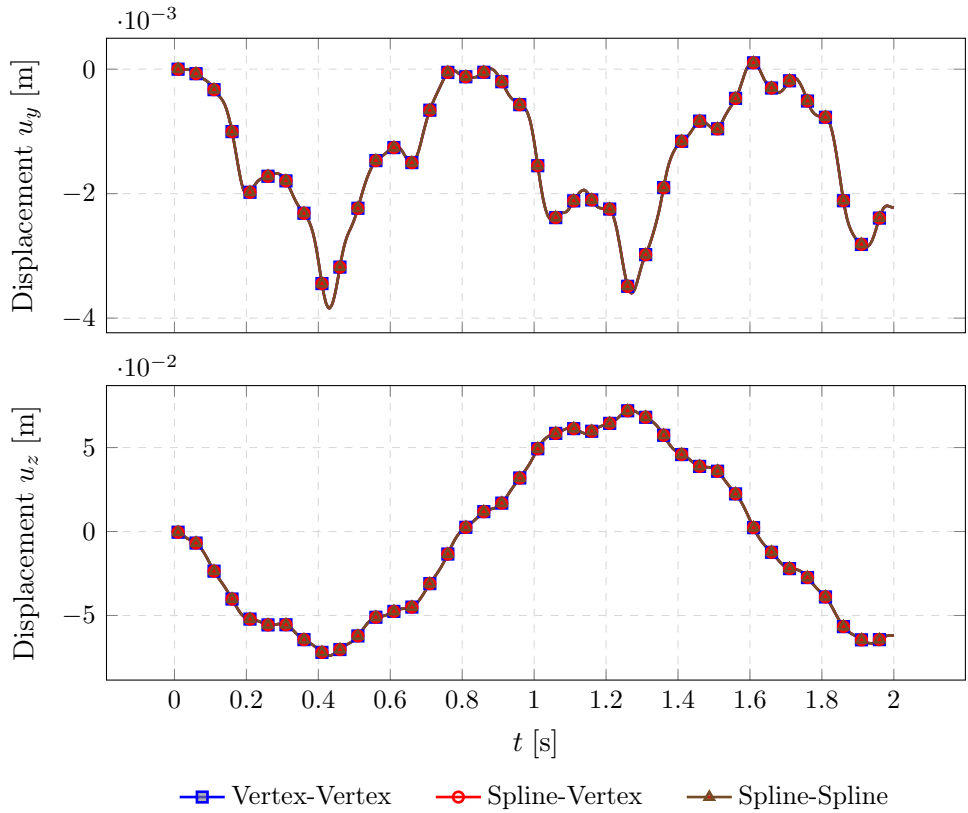

Figure~\ref{figure:nonlinear-vertical-beam-tip-displacement} shows the tip displacement over time under a nonlinear load. The spline-spline based, spline-vertex based, and vertex-vertex based coupling methods yield consistent results under complex loading conditions.
\FloatBarrier

\section{Communication Overhead Comparison and Analysis}
\label{section:communication-overhead}
In partitioned multiphysics simulations, communication overhead is a critical metric for evaluating coupling method efficiency. We conducted both theoretical analysis and experimental validation of communication overhead for different coupling methods to assess the advantages of spline-based coupling methods compared to traditional vertex-based approaches.

\subsection{Theoretical Communication Overhead}

\begin{figure}[htbp]
\centering
\begin{minipage}[b]{0.48\linewidth}
    \vspace*{0pt}
    \centering
    \begin{tikzpicture}[]
    \begin{axis}[
        title={$N_t=10000$},
        xlabel={\# subdivisions $r$ per direction},
        ylabel={Total Communication Overhead (KB)},
        xmode=log,
        xtick={2,4,8,16,32,64,128},
        xticklabels={$2^1$,$2^2$,$2^3$,$2^4$,$2^5$,$2^6$,$2^7$},
        width=\linewidth,
        height=8cm,
        tick style={black},
        cycle list name=customcolorlist,
        legend cell align={left}
    ]
    \pgfplotsinvokeforeach{2,3,4,5}{%
        \addplot+[only marks, domain=2:128, samples=7] 
            {(10000*((x+#1)^2*3) + 2*(x + 2*#1 + 1))*8/1024};
        \addlegendentry{\(O_{\text{spline}},\; p=#1\)}

        \addplot+[no markers, domain=2:128, samples=7] 
            {(10000*(x^2*(#1+1)^2))*8/1024};
        \addlegendentry{\(O_{\text{vertex}},\; p=#1\)}
    }
    \end{axis}
    \end{tikzpicture}
    \label{fig:big}
\end{minipage}%
\hfill
\begin{minipage}[b]{0.45\linewidth}
    \centering
        \centering
        \begin{tikzpicture}[]
        \begin{groupplot}[
            group style={group size=1 by 2,
            x descriptions at=edge bottom,
            y descriptions at=edge left,
            horizontal sep=0.5cm,
            vertical sep=1cm,},
            xlabel={\# subdivisions $r$ per direction},
            ylabel={Overhead (KB)},
            xmode=log,
            xtick={2,4,8,16,32,64,128},
            xticklabels={$2^1$,$2^2$,$2^3$,$2^4$,$2^5$,$2^6$,$2^7$},
            width=\linewidth,
            height={4.3cm},
            tick style={black},
            every axis plot/.append style={line width=1pt},
            cycle list name=customcolorlist,
        ]
         \nextgroupplot[]
         \pgfplotsinvokeforeach{2,3,4,5}{%
            \addplot+[only marks, domain=2:128, samples=7] 
                {(((x + #1)^2*3) + 2*(x + 2*#1 + 1))*8/1024};
            \addplot+[no markers, domain=2:128, samples=7] 
                {(x^2*(#1+1)^2)*8/1024};
        }
        \node[fill=white] at (rel axis cs: 0.2,0.8) {$N_t=1$};
        \nextgroupplot[]
        \pgfplotsinvokeforeach{2,3,4,5}{%
            \addplot+[only marks, domain=2:128, samples=7] 
                {(100*((x+ #1)^2*3) + 2*(x + 2*#1 + 1))*8/1024};
            \addplot+[no markers, domain=2:128, samples=7] 
                {(100*(x^2*(#1+1)^2))*8/1024};
        }
        \node[fill=white] at (rel axis cs: 0.2,0.8) {$N_t=100$};
        \end{groupplot}
        \end{tikzpicture}
        \label{fig:small1}
\end{minipage}
\caption{Comparison of the communication overhead between vertex-based (FEM) and spline-based (IGA) coupling. The horizontal axis represents the number of parameter subdivisions per direction, serving as an equivalent indicator of the unique knots in IGA while reflecting mesh resolution in FEM. The vertical axis shows the total communication overhead in kilobytes.}
\label{fig:theoretical-communication-combined}
\end{figure}
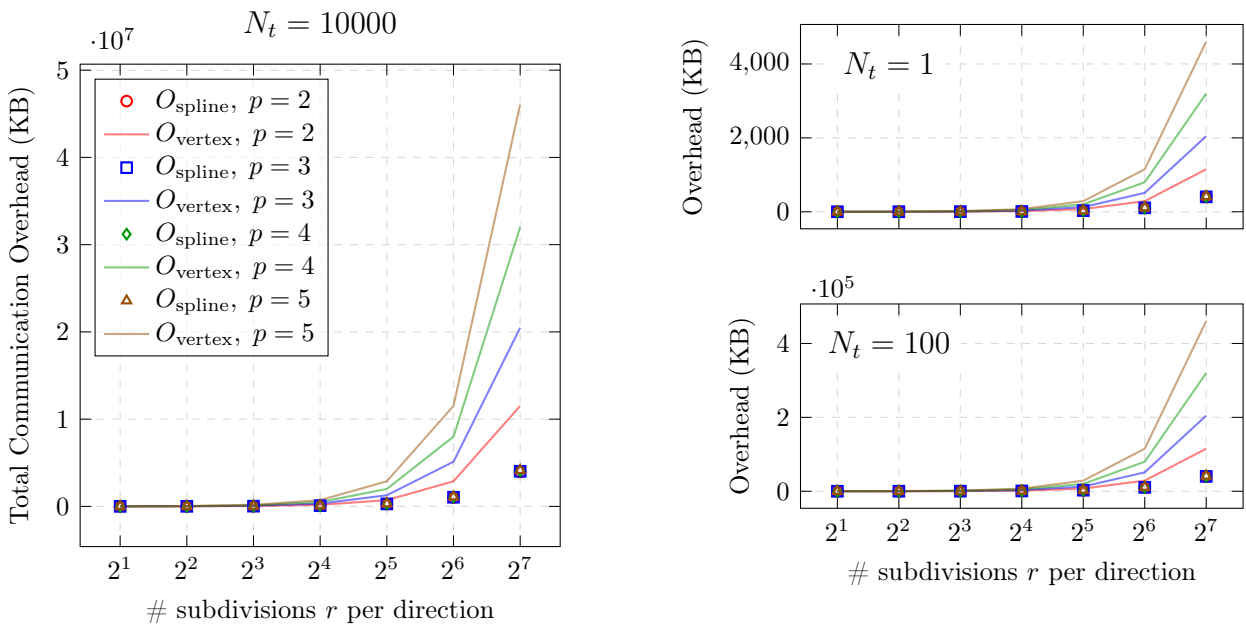

As described in Section~\ref{subsec:spline-based-coupling}, the communication overhead for spline-based coupling methods primarily depends on the number of control points and parametric representation, while vertex-based methods see overhead increase exponentially with mesh refinement and quadrature point quantities. Figure~\ref{fig:theoretical-communication-combined} illustrates this theoretical difference.

For vertex-based coupling using traditional Gauss quadrature with $p+1$ quadrature points per parametric direction per element, the communication overhead $O$ is proportional to the number of interface elements $n_{\text{element}}$, the polynomial degree $p$, the number of timesteps $N_t$, and the parametric dimension $\hat{d}$:
\begin{equation}
\label{eq:overhead_vertex}
O_{\text{quad}} = N_t \cdot d \cdot n_{\text{element}} \cdot (p+1)^{\hat{d}} 
\end{equation}

For spline-based coupling, the total communication overhead is primarily governed by the number of control points, the physical field dimension $d$, and the polynomial degrees $p_i$ along each parametric direction. Unlike mesh-based coupling (which requires per-node communication), spline-based coupling grows at a much slower rate:
\begin{equation}
O_{\text{spline}} = N_t \cdot d \cdot \left( \prod_{i=1}^{\hat{d}} n_i \right) + \sum_{i=1}^{\hat{d}} (n_i + p_i + 1) \label{eq:overhead_spline}
\end{equation}

This expression consists of two main components:
\begin{itemize}
\item $N_t \cdot d \cdot \left( \prod_{i=1}^{\hat{d}} n_i \right)$: the total number of scalar values exchanged during time integration. Here, $n_i$ denotes the number of control points along the $i$-th parametric direction (with $\hat{d}$ such directions in total), and $d$ represents the number of coupled physical fields.

\item $\sum_{i=1}^{\hat{d}} (n_i + p_i + 1)$: the cumulative size of the knot vectors in each parametric direction. This geometric data is exchanged once during initialization and remains unchanged throughout the simulation.
\end{itemize}

We use preCICE \cite{Chourdakis2021} to mediate the exchange of field data (e.g., displacement, velocity, pressure) between solvers. On the interface, both the geometry and field values are represented using B-splines or NURBS. The spline representation allows us to transmit compact and structured data, where:
\begin{itemize}
    \item The knot vector (of length $\sum_{i=1}^d (n_{i} + p_i + 1)$) is transferred at the beginning of the simulation.
    \item The field values are exchanged at every coupling step through the control points.
\end{itemize}

This comparison highlights a major advantage of spline-based coupling. As shown in Figure~\ref{fig:theoretical-communication-combined}, vertex-based methods lead to an exponential increase in communication overhead as mesh refinement and polynomial degree grow, making them computationally expensive for high-fidelity simulations. In contrast, spline-based methods maintain significantly lower communication costs, as their overhead scales more efficiently with refinement.

\FloatBarrier

    

\FloatBarrier

The analysis shows that as the mesh becomes more refined, the difference in communication overhead between vertex-based and spline-based methods will grow significantly. 

\subsection{Experimental Validation Results}
The experimental results confirm our theoretical predictions about communication overhead in different coupling methods.
We measured both the actual data transfer volume during the initial timestep ($N_t=1$) and total communication time per dimension over $N_t=10000$ timestep in a vertical beam test case across various mesh refinement levels. 

Our instrumentation specifically targeted only the coupling-related data transfer functions, ensuring that the measured values reflect purely the communication overhead of different coupling methods without contamination from other computational processes. 

In the case of spline-based coupling, additional overhead arises from transmitting knot vectors alongside control point data. To support tensor-product B-spline representations in multiple parametric directions, knot vectors are embedded into a unified matrix structure with \texttt{NaN} padding. Since the knot vectors are defined separately for each parametric direction, we embed them into a unified matrix where each row corresponds to one direction. \texttt{NaN} values are used to fill the remaining entries. This ensures consistency in data formatting and facilitates reliable transfer and storage. 

Here we consider a parametric dimension $\hat{d} = 2$ and a physical dimension $d = 3$. The tensor-product B-spline space is defined by two univariate knot vectors associated with each parametric direction:

$$
\Xi = \{\xi_1, \xi_2, \ldots, \xi_{n_\xi + p_\xi + 1}\} \subset \mathbb{R}, \quad
\mathscr{H} = \{\eta_1, \eta_2, \ldots, \eta_{n_\eta + p_\eta + 1}\} \subset \mathbb{R},
$$
where each knot vector is a non-decreasing sequence of real numbers. Here, $n_\xi$ and $n_\eta$ denote the number of basis functions in the $\xi$- and $\eta$-directions, respectively, and $p_\xi$, $p_\eta$ are the corresponding polynomial degrees.

In the implementation, while arranging knot vectors into a unified matrix format to ensure consistent data transmission, we inevitably introduce some overhead.

Since the lengths of knot vectors can differ between parametric directions, \texttt{NaN} values are used as padding to maintain consistent matrix dimensions:
$$
\texttt{knotMatrix} =
\begin{bmatrix}
\xi_1 & \xi_2 & \dots & \xi_{n_\xi + p_\xi +1} & \texttt{NaN} & \dots & \texttt{NaN} \\
\texttt{NaN} & \dots & \texttt{NaN} & \eta_1 & \eta_2 & \dots & \eta_{n_\eta + p_\eta + 1}
\end{bmatrix}
$$
In the numerical experiment, we consider the symmetric case where $n_\xi = n_\eta = n$ and $p_\xi = p_\eta = p$, so that both knot vectors $\boldsymbol{\Xi}$ and $\mathscr{H}$ have the same length, namely $n + p + 1$. To embed them into a unified matrix for transfer, we allocate a matrix of size $2 \times (2n + 2p + 2)$, where each row stores one knot vector padded with \texttt{NaN} values in the complementary half. Consequently, the total number of \texttt{NaN} entries equals $2(n + p + 1)$.

All timing data were collected on a single compute node of the DelftBlue supercomputer \cite{DHPC2024}. Given that runtime performance can be affected by external processes and background system noise, all measurements were repeated ten times in an isolated environment with minimal external interference. The values reported represent the averaged results, ensuring reliability and minimizing the influence of transient fluctuations.

\FloatBarrier
\begin{figure}[h!]
    \centering
    \begin{subfigure}[t]{0.45\textwidth}
        \centering
        \begin{tikzpicture}[]
        \begin{axis}[
            width=\linewidth,
            height=8cm,
            xmode=log,
            log basis x=2,
            xlabel={\# subdivisions $r$ per direction},
            ylabel={Vertex (quadarture) communication size (KB)},
            cycle list name=customcolorlist
        ]
        
        \addplot+[only marks]
        coordinates {(2,0.84375)(4,3.375)(8,13.5)(16,54)(32,216)(64,864)(128,3456)};
        \addlegendentry{Exp. $p=2$}

        \addplot+[no markers]
        coordinates {(2,0.84375)(4,3.375)(8,13.5)(16,54)(32,216)(64,864)(128,3456)};
        \addlegendentry{An. $p=2$}
        
        \addplot+[only marks]
        coordinates {(2,1.5)(4,6)(8,24)(16,96)(32,384)(64,1536)(128,6144)};
        \addlegendentry{Exp. $p=3$}

        \addplot+[no markers]
        coordinates {(2,1.5)(4,6)(8,24)(16,96)(32,384)(64,1536)(128,6144)};
        \addlegendentry{An. $p=3$}
        
        \addplot+[only marks]
        coordinates {(2,2.34375)(4,9.375)(8,37.5)(16,150)(32,600)(64,2400)(128,9600)};
        \addlegendentry{Exp. $p=4$}

        \addplot+[no markers]
        coordinates {(2,2.34375)(4,9.375)(8,37.5)(16,150)(32,600)(64,2400)(128,9600)};
        \addlegendentry{An. $p=4$}
        
        \addplot+[only marks]
        coordinates {(2,3.375)(4,13.5)(8,54)(16,216)(32,864)(64,3456)(128,13824)};
        \addlegendentry{Exp. $p=5$}

        \addplot+[no markers]
        coordinates {(2,3.375)(4,13.5)(8,54)(16,216)(32,864)(64,3456)(128,13824)};
        \addlegendentry{An. $p=5$}
        \end{axis}
        \end{tikzpicture}
        \caption{Vertex-based communication overhead $N_t=1$.}
        \label{fig:comm-overhead-exp-theo-vertex}
    \end{subfigure}
    \hfill
    \begin{subfigure}[t]{0.45\textwidth}
        \centering
        \begin{tikzpicture}[]
        \begin{axis}[
            width=\linewidth,
            height=8cm,
            xmode=log,
            log basis x=2,
            xlabel={\# subdivisions $r$ per direction},
            ylabel={Spline communication size (KB)},
            cycle list name=customcolorlist
        ]
        
        \addplot+[only marks]
        coordinates {(2, 0.59375) (4, 1.125) (8, 2.75) (16, 8.25) (32, 28.25) (64, 104.25) (128, 400.25)};
        \addlegendentry{Exp. $p=2$}

        \addplot+[no markers, domain=2:128, samples=7]
            {(((x+ 2)^2*3) + 2*(x + 2*2 + 1))*8/1024};
        \addlegendentry{An. $p=2$}
        
        \addplot+[only marks]
        coordinates {(2, 0.867188) (4, 1.49219) (8, 3.30469) (16, 9.17969) (32, 29.9297) (64, 107.43) (128, 406.43)};
        \addlegendentry{Exp. $p=3$}

        \addplot+[no markers, domain=2:128, samples=7]
            {(((x+ 3)^2*3) + 2*(x + 2*3 + 1))*8/1024};
        \addlegendentry{An. $p=3$}
        
        \addplot+[only marks]
        coordinates {(2, 1.1875) (4, 1.90625) (8, 3.90625) (16, 10.1562) (32, 31.6562) (64, 110.656) (128, 412.656)};
        \addlegendentry{Exp. $p=4$}

        \addplot+[no markers, domain=2:128, samples=7]
            {(((x+ 4)^2*3) + 2*(x + 2*4 + 1))*8/1024};
        \addlegendentry{An. $p=4$}

        \addplot+[only marks]
        coordinates {(2, 1.55469) (4, 2.36719) (8, 4.55469) (16, 11.1797) (32, 33.4297) (64, 113.93) (128, 418.93)};
        \addlegendentry{Exp. $p=5$}

        \addplot+[no markers, domain=2:128, samples=7]
            {(((x+ 5)^2*3) + 2*(x + 2*5 + 1))*8/1024};
        \addlegendentry{An. $p=5$}

        \end{axis}
        \end{tikzpicture}
        \caption{Spline-based communication overhead $N_t=1$.}
        \label{fig:comm-overhead-exp-theo-spline}
    \end{subfigure}

    \vspace{0.5em}

    \hfill 
    \begin{subfigure}[t]{\linewidth}
        \centering
        \begin{tikzpicture}[]
        \begin{axis}[
            ybar,                           
            bar width=5pt,                  
            enlargelimits=0.15,
            symbolic x coords={2,4,8,16,32,64,128},
            xtick=data,
            xticklabels={$2^1$,$2^2$,$2^3$,$2^4$,$2^5$,$2^6$,$2^7$},
            ylabel={Dimensionless Communication Time},
            xlabel={\# subdivisions $r$ per direction},
            width=\linewidth,
            height=8cm,
            cycle list name=customcolorlist,
            every axis plot/.append style={fill,draw=none,no markers},
            legend style={
                at={(0.5,-0.25)},
                anchor=north,
                legend columns=4,
                column sep=10pt,
                font=\footnotesize
            },
            no markers
        ]
            \addplot+[bar shift=-15pt] coordinates {
                (2, 1.000)
                (4, 1.1486)
                (8, 1.1886)
                (16, 1.3214)
                (32, 1.0016)
                (64, 1.0801)
                (128, 1.1574)
            };
            \addlegendentry{Spline $p=2$}
            
                \addplot+[bar shift=9pt] coordinates {
                (2, 1.000)
                (4, 1.2062)
                (8, 1.6971)
                (16, 2.6587)
                (32, 4.7173)
                (64, 9.2688)
                (128, 16.8319)
            };
            \addlegendentry{Vertex $p=2$}
            \addplot+[bar shift=-9pt] coordinates {
                (2, 1.000)
                (4, 0.7920)
                (8, 0.8887)
                (16, 0.9484)
                (32, 0.9387)
                (64, 1.0570)
                (128, 1.1574)
            };
            \addlegendentry{Spline $p=3$}
                \addplot+[bar shift=15pt] coordinates {
                (2, 1.000)
                (4, 1.3416)
                (8, 2.4116)
                (16, 4.4658)
                (32, 8.5610)
                (64, 18.4161)
                (128, 34.9442)
            };
            \addlegendentry{Vertex $p=3$}
            
            \addplot+[bar shift=-3pt] coordinates {
                (2, 1.000)
                (4, 0.7804)
                (8, 0.9711)
                (16, 0.9901)
                (32, 0.9359)
                (64, 1.0705)
                (128, 1.1677)
            };
            \addlegendentry{Spline $p=4$}
            
            \addplot+[bar shift=21pt] coordinates {
                (2, 1.000)
                (4, 1.5494)
                (8, 2.6280)
                (16, 4.9168)
                (32, 10.0070)
                (64, 20.2998)
                (128, 34.9133)
            };
            \addlegendentry{Vertex $p=4$}
            
            \addplot+[bar shift=3pt] coordinates {
                (2, 1.000)
                (4, 1.0378)
                (8, 1.0131)
                (16, 1.0185)
                (32, 1.0239)
                (64, 1.1112)
                (128, 1.2115)
            };
            \addlegendentry{Spline $p=5$}
            

            \addplot+[bar shift=27pt] coordinates {
                (2, 1.000)
                (4, 1.4788)
                (8, 2.8776)
                (16, 5.4570)
                (32, 11.7660)
                (64, 28.4937)
                (128, 56.5501)
            };
            \addlegendentry{Vertex $p=5$}
            
            \end{axis}
            \end{tikzpicture}
            \caption{Communication Time (normalized by $r{=}2^1$).}
        \label{fig:comm-overhead-exp-time}
    \end{subfigure}
    \caption{Experimental (Exp.) and theoretical (An.) comparison of communication overhead with respect to the number of parameter subdivisions $r$ per direction, with for $r\in\{1,2,...,7\}$
    .}
    \label{fig:communication-overhead-comparison}
\end{figure}
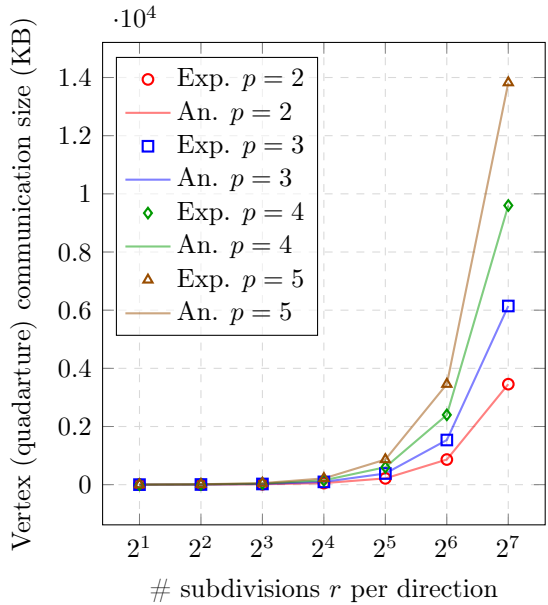
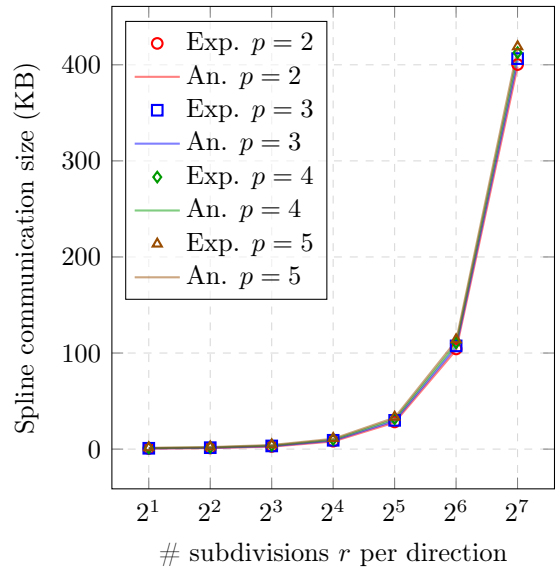
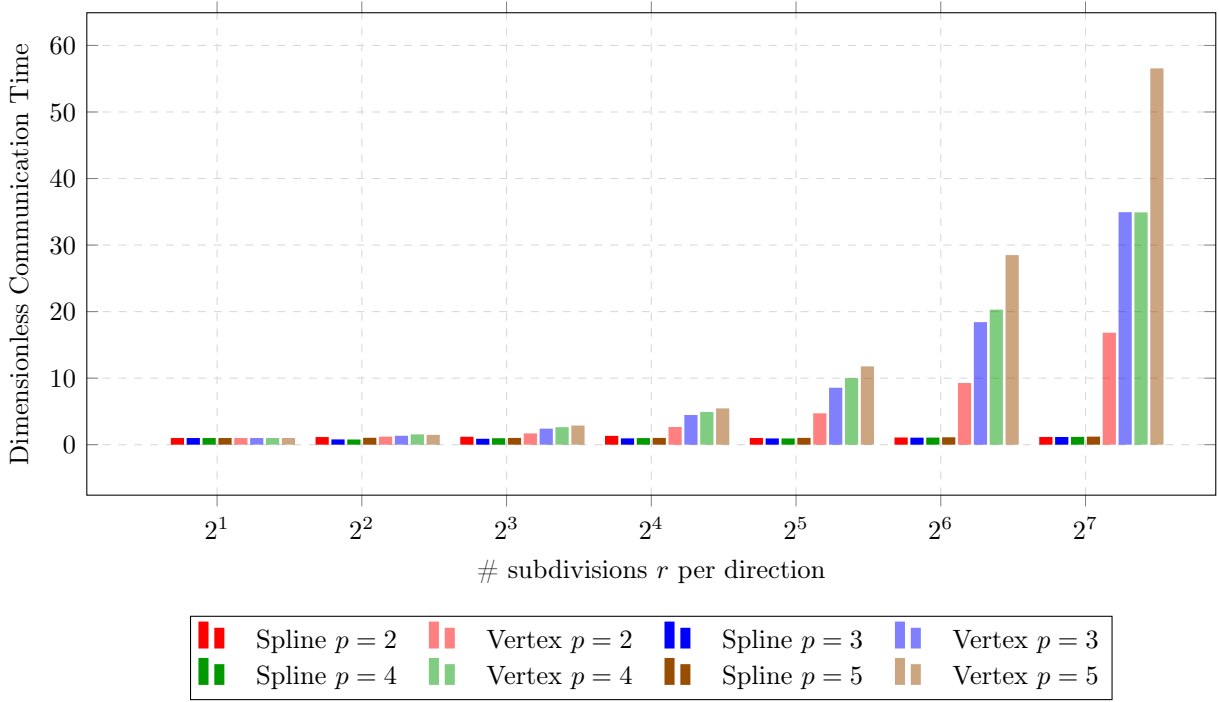
\FloatBarrier
Figure~\ref{fig:communication-overhead-comparison} presents a comprehensive analysis of communication overhead for different coupling methods through four complementary perspectives. 
In Figure~\ref{fig:comm-overhead-exp-theo-vertex}, we compare experimental measurements against theoretical predictions for vertex-based coupling in terms of data volume (KB), demonstrating close agreement between the predicted quadratic growth and actual performance across increasing mesh refinement levels.

Figure~\ref{fig:comm-overhead-exp-theo-spline} provides a similar data volume comparison for spline-based coupling, where both theoretical and experimental results confirm the significantly reduced growth rate in communication overhead as mesh resolution increases. 
Finally, Figure~\ref{fig:comm-overhead-exp-time} translates these data volume differences into actual communication time measurements, validating that the theoretical efficiency gains of spline-based coupling translate to meaningful performance improvements in practice, with timing differences becoming increasingly pronounced at finer mesh resolutions.

As shown in Figure~\ref{fig:comm-overhead-exp-theo-spline}, small discrepancies are observed between the experimental and theoretical values for spline-based coupling. These deviations are primarily attributed to the inclusion of \texttt{NaN}-padded entries when transmitting knot vectors, which slightly increased the actual communication size beyond the idealized theoretical estimates. 

Table~\ref{tab:exact_error} supports this observation by listing the absolute differences between experimental and analytical spline-based communication overheads at $r=7$ for various polynomial degrees $p$, corresponding to the results shown in Figure~\ref{fig:comm-overhead-exp-theo-spline}. The absolute discrepancy is computed as
\begin{equation}
    \varepsilon_{\text{absolute}} = \left| O_{\text{exp}} - O_{\text{an}} \right|,
\end{equation}
where $O_{\text{exp}}$ denotes the measured (experimental) communication overhead in kilobytes, and $O_{\text{an}}$ is the corresponding analytical estimate based on Equation~\ref{eq:overhead_spline}.

\begin{table}[h!]
\centering
\caption{Absolute discrepancy (KB) between measured and analytical communication sizes at $r=7$, with corresponding analytical number of \texttt{NaN}-padded entries.}
\label{tab:exact_error}
\begin{tabular}{c|c|c|c}
\toprule
$p$ & Absolute Discrepancy (KB) & \# Measured \texttt{NaN} entries & \# Analytical \texttt{NaN} entries\textsuperscript{*}\tablefootnote{Computed as $2(n + p + 1)$, where $n = 2^r$ is the number of basis functions per parametric direction and $p$ is the spline degree.} \\
\midrule
2 & 2.078125 & 266.00 & 266 \\
3 & 2.1096875 & 270.04 & 270 \\
4 & 2.140375 & 273.97 & 274 \\
5 & 2.1721875 & 278.04 & 278 \\
\bottomrule
\end{tabular}
\vspace{0.5em}
\begin{flushleft}
\end{flushleft}
\end{table}

\FloatBarrier

Key observations from comparing the experimental results (Figure~\ref{fig:communication-overhead-comparison}) with theoretical results (Figure~\ref{fig:theoretical-communication-combined}) include:

\begin{enumerate}
    \item The communication overhead for Vertex-based coupling methods increases rapidly with mesh refinement in practice, which is consistent with the theoretical prediction.
    \item Spline-based coupling methods (including spline-vertex and IGA-IGA approaches) demonstrate significantly lower growth rates in communication overhead, confirming the theoretical advantages.
    \item The experimental curves closely match the theoretical predictions, although the absolute values differ due to implementation details and overhead from the coupling library.
\end{enumerate}

These results confirm the communication efficiency advantages of spline-based coupling methods, particularly when dealing with high-resolution interface meshes. The close agreement between theoretical predictions and experimental measurements validates our communication overhead model and provides a solid foundation for estimating performance in larger, more complex simulations.


\section{Benchmark problems}
\label{section:results}
 So far we have demonstrated the superior communication efficiency of spline-based coupling methods in solvers that support spline representation, while also illustrating how vertex-based coupling methods can effectively interface with solvers utilizing alternative discretization techniques through the preCICE library. In this section, we extend our analysis to validate the accuracy and robustness of the proposed methods through a series of carefully selected benchmark problems that progressively increase in complexity and physical relevance. The solid domain is modeled using the G+Smo library \cite{gismo_paper}, which provides native support for IGA, while the fluid domain is solved using OpenFOAM\cite{OpenFOAM}, a widely used finite volume-based CFD framework. 
 
\subsection{Partitioned Heat Conduction}
\FloatBarrier
In this example, we employed the spline-based communication method to simulate heat transfer across partitioned domains. 
Beyond the previously demonstrated advantages in communication efficiency and information preservation, this application reveals two additional benefits particularly relevant to FSI. Let's first start with the benchmark setup. 

\begin{figure}[h!]
    \centering
    \includegraphics[width=0.8\linewidth]{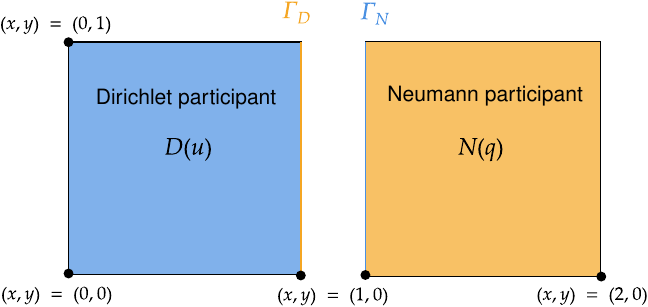}
    \caption{Partitioned heat conduction problem setup.}
    \label{fig:partitioned-heat-conduction-setup}
\end{figure}

The partitioned heat conduction setup, shown in Figure~\ref{fig:partitioned-heat-conduction-setup}, consists of two square regions with different thermal properties. The governing equation is the dimensionless, transient heat equation in two spatial dimensions:
\begin{equation} \label{eq:heat-equation}
\begin{aligned}
    \frac{\partial u}{\partial t} &= \nabla^2 u + f && \text{in } \Omega = [0,2] \times [0,1],\quad t \in [0,T], \\
    u &= g && \text{on } \partial \Omega,\quad t \in [0,T],
\end{aligned}
\end{equation}
where \( u(x, y, t) \) is the temperature field, \( f \) is a source term, and \( g \) denotes prescribed boundary data.

To validate the numerical method, we adopt a manufactured solution of the form
\begin{equation} \label{eq:exact-solution}
    u_{\text{exact}}(x, y, t) = 1 + x^2 + \alpha y^2 + \beta t,
\end{equation}
with constants \( \alpha, \beta \in \mathbb{R} \). Substituting \( u_{\text{exact}} \) into \eqref{eq:heat-equation} yields the corresponding source term \( f \) and Dirichlet boundary condition \( g \).

Following the formulation in \cite{heat_conduction}, the domain \( \Omega \) is partitioned into two non-overlapping subdomains:
\[
\Omega_D = [0,1] \times [0,1], \qquad \Omega_N = [1,2] \times [0,1],
\]
such that \( \Omega = \Omega_D \cup \Omega_N \). These subdomains are coupled along the shared interface
\[
\Gamma = \Omega_D \cap \Omega_N = \{1\} \times [0,1],
\]
as shown in Figure~\ref{fig:partitioned-heat-conduction-setup}. This configuration is referred to as the partitioned heat equation.

During the coupling procedure, either the temperature \( u \) or the normal component of the heat flux
\[
q_n := \nabla u \cdot \mathbf{n}
\]
is exchanged across \( \Gamma \). Specifically, one subdomain applies a Dirichlet condition
\[
u = u_\Gamma \quad \text{on } \Gamma,
\]
while the other enforces a Neumann condition
\[
\nabla u \cdot \mathbf{n} = q_\Gamma \quad \text{on } \Gamma,
\]
where \( u_\Gamma \) and \( q_\Gamma \) represent the interface temperature and heat flux, respectively.

\FloatBarrier

\textbf{Naturally non-conforming support:} One of the main advantages of the proposed coupling strategy is its natural ability to handle non-conforming discretizations at the coupling interface. Even when the mesh is geometrically discontinuous across the interface, the use of spline-based representations allows us to preserve functional smoothness. In the current implementation, non-conforming meshes are assumed to be nested, meaning that the control points or knot vectors of one mesh are fully contained within those of the other. 


\FloatBarrier
\begin{figure}[htbp]
    \centering
    \begin{subfigure}[t]{0.45\linewidth}
        \centering
        \includegraphics[width=\linewidth,trim={120 150 120 175},clip]{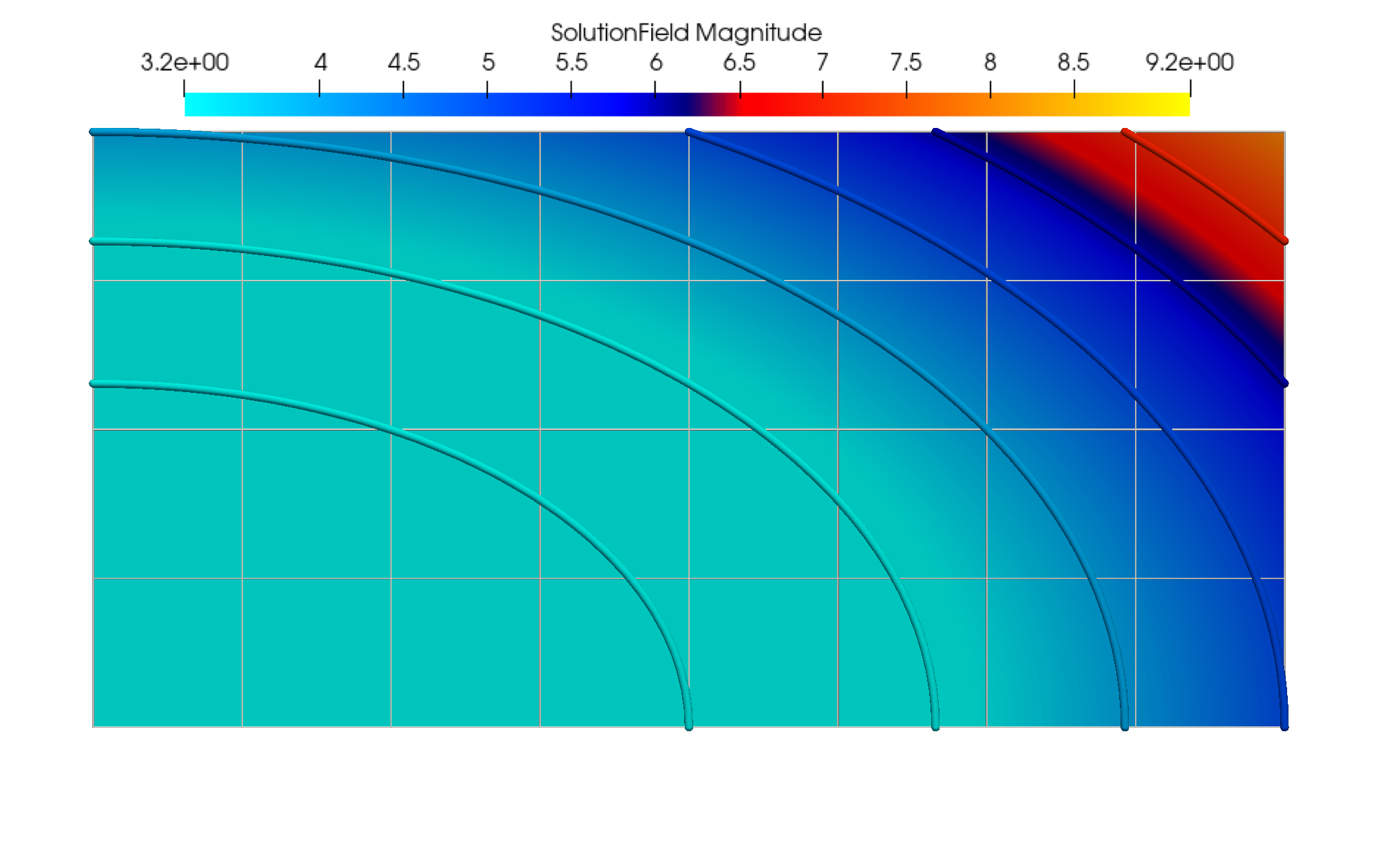}
        \caption{$t=0\,\mathrm{s},\;r_{N}=2$}
        \label{fig:partitioned-heat-gridDa}
    \end{subfigure}
    \hfill
    \begin{subfigure}[t]{0.45\linewidth}
        \centering
        \includegraphics[width=\linewidth,trim={120 150 120 175},clip]{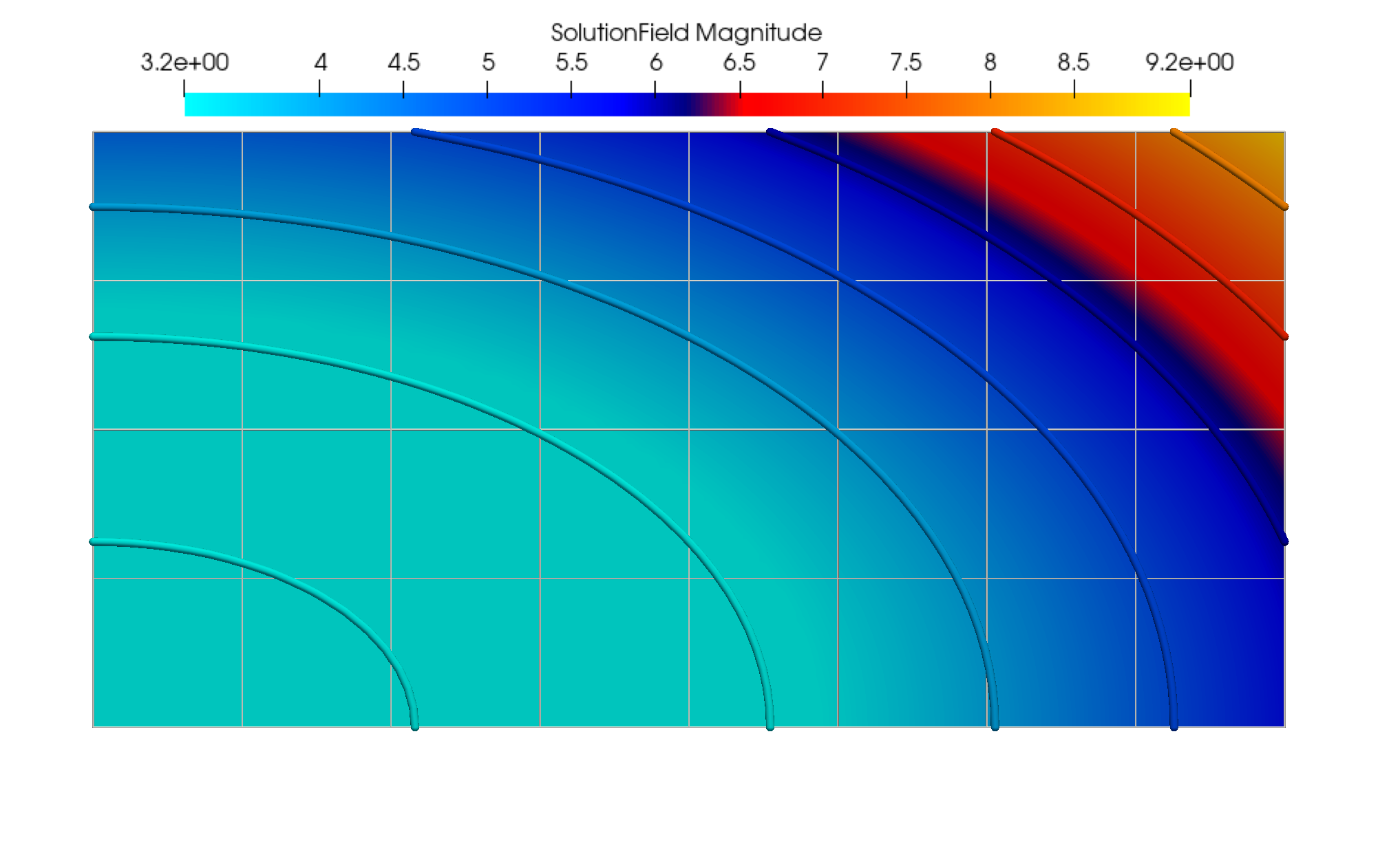}
        \caption{$t=0.6\,\mathrm{s},\;r_{N}=2$}
    \end{subfigure}

    \begin{subfigure}[t]{0.45\linewidth}
        \centering
        \includegraphics[width=\linewidth,trim={120 150 120 175},clip]{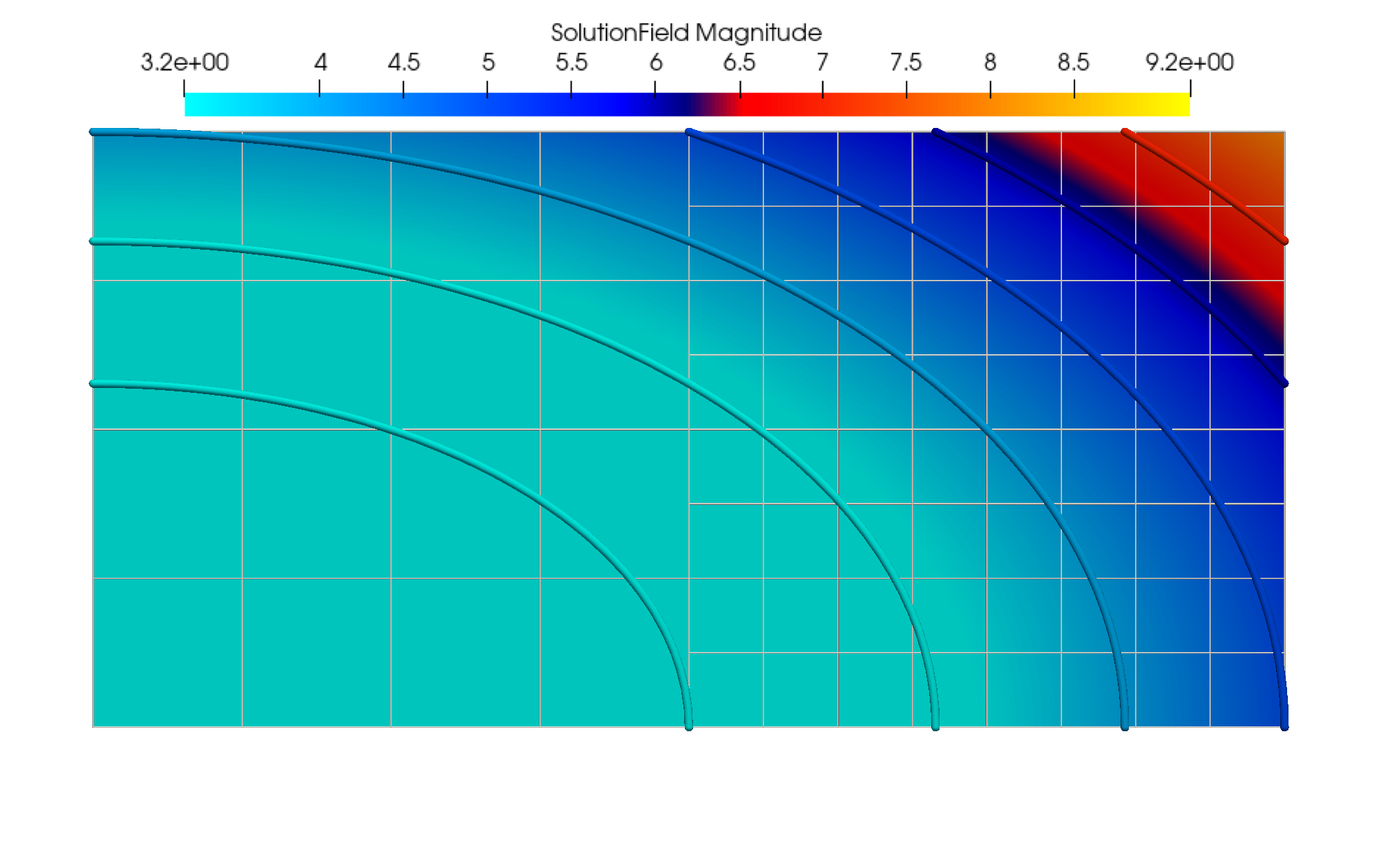}
        \caption{$t=0\,\mathrm{s},\;r_{N}=3$}
        \label{fig:partitioned-heat-griddb}
    \end{subfigure}
    \hfill
    \begin{subfigure}[t]{0.45\linewidth}
        \centering
        \includegraphics[width=\linewidth,trim={120 150 120 175},clip]{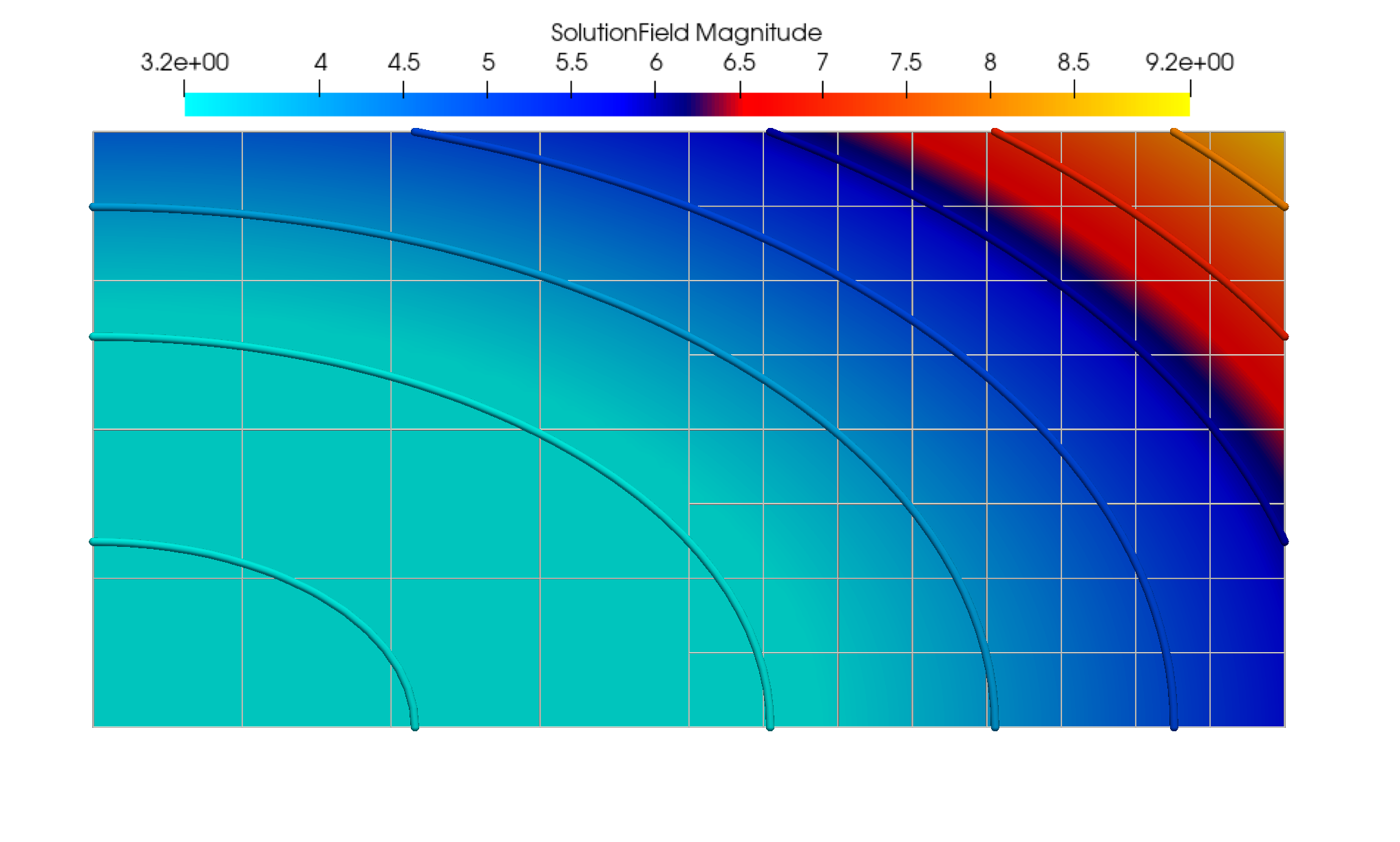}
        \caption{$t=0.6\,\mathrm{s},\;r_{N}=3$}
    \end{subfigure}

    \begin{subfigure}[t]{0.45\linewidth}
        \centering
        \includegraphics[width=\linewidth,trim={120 150 120 175},clip]{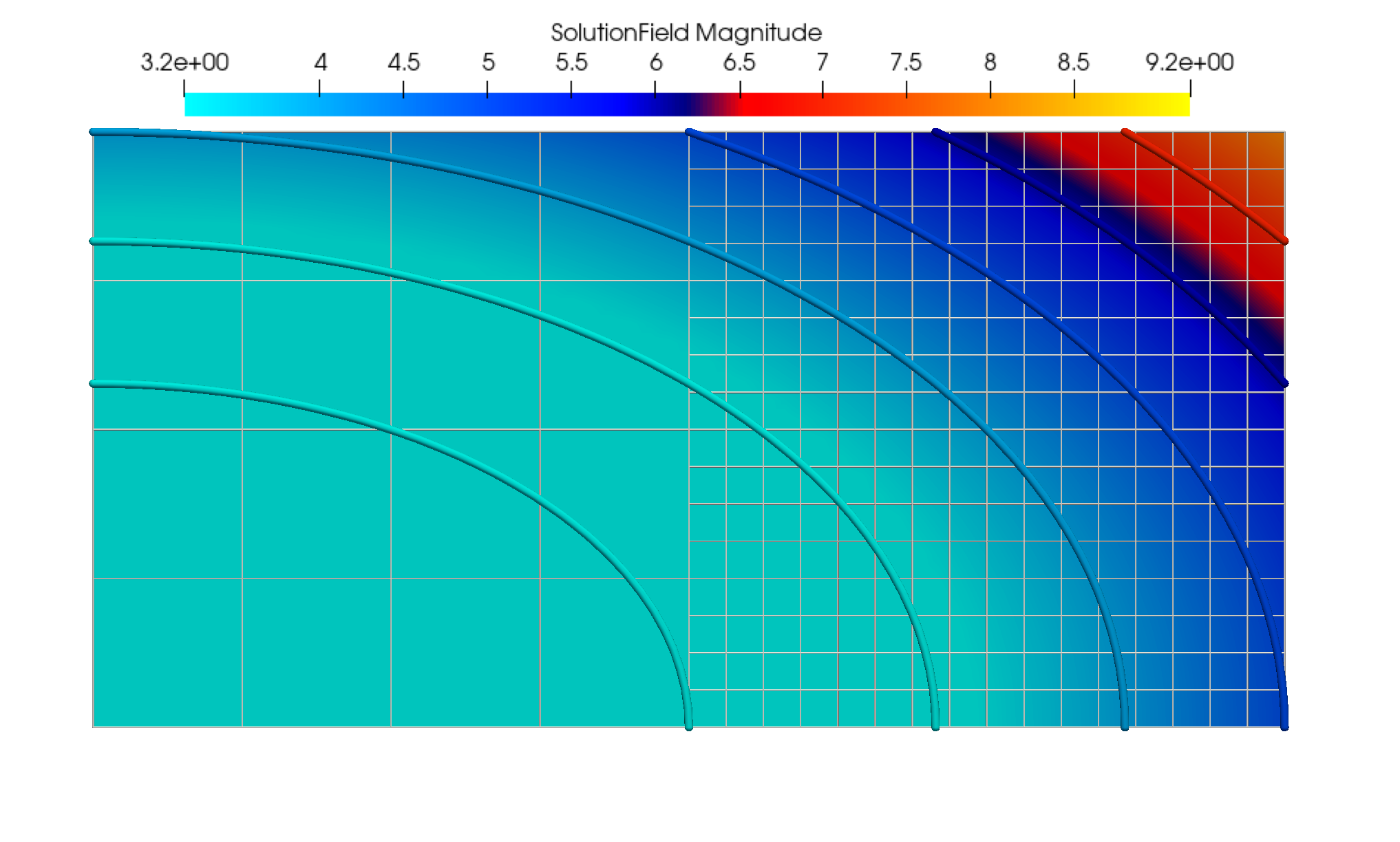}
        \caption{$t=0\,\mathrm{s},\;r_{N}=4$}
        \label{fig:partitioned-heat-griddc}
    \end{subfigure}
    \hfill
    \begin{subfigure}[t]{0.45\linewidth}
        \centering
        \includegraphics[width=\linewidth,trim={120 150 120 175},clip]{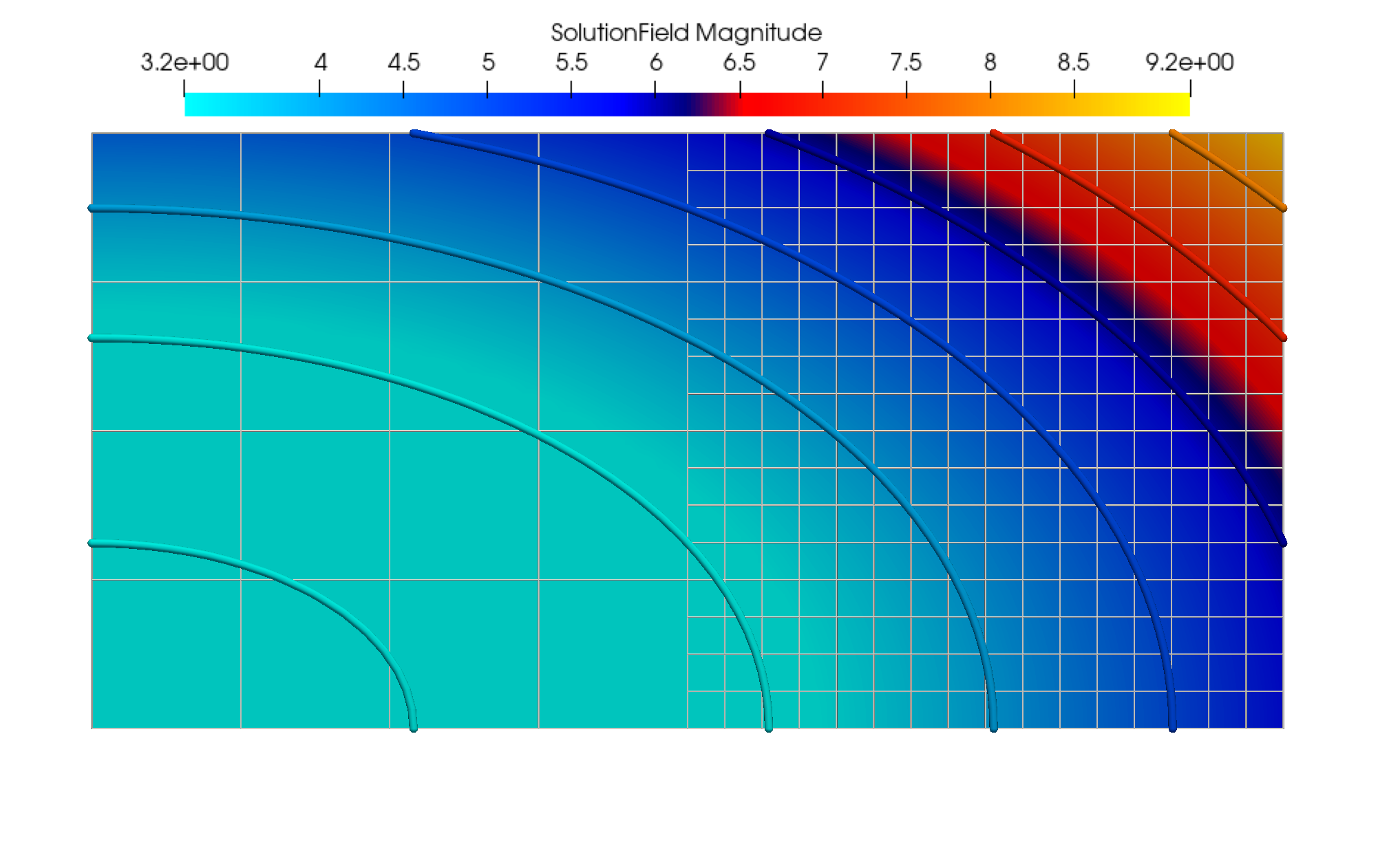}
        \caption{$t=0.6\,\mathrm{s},\;r_{N}=4$}
    \end{subfigure}

    \begin{subfigure}[t]{\linewidth}
        \centering
        \includegraphics[width=\linewidth,trim={0 1050 0 0},clip]{figures/partitioned-heat0r2.png}
    \end{subfigure}

    \caption{Temperature field of a partitioned spline-mesh domain with varying $\Omega_N$ refinement levels \(r_{N}\).}
    \label{fig:partitioned-heat-grid}
\end{figure}
\FloatBarrier

As illustrated in Figure~\ref{fig:partitioned-heat-grid}, the Dirichlet subdomain $\Omega_D$ is discretized with a coarser mesh refined once in Figure~\ref{fig:partitioned-heat-gridDa} and twice in Figure~\ref{fig:partitioned-heat-griddb} and Figure~\ref{fig:partitioned-heat-griddc}, while the Neumann subdomain $\Omega_{N}$ employs progressively finer meshes. Since the coupling interface is represented as a smooth spline curve, physical quantities such as temperature or heat flux can be consistently evaluated at any point on the interface. This enables direct and accurate data exchange between nonmatching discretizations, without requiring interpolation or projection, and ensures continuity of the thermal field across the interface.

\textbf{Inherited derivatives via spline representation:} One of the distinctive advantages of using spline-based representations for coupling is that the derivatives of the solution field are inherently preserved and transferred across the interface. Since $\Gamma_D$ is represented as a spline function, its smoothness and continuity naturally extend to its derivatives. As a result, when exchanging data across non-conforming interfaces, both the field values and their gradients are consistently evaluated from the same spline representation, without the need for additional reconstruction or projection.
\FloatBarrier

\begin{figure}[h!]
    \centering
    \includegraphics[width=0.9\linewidth]{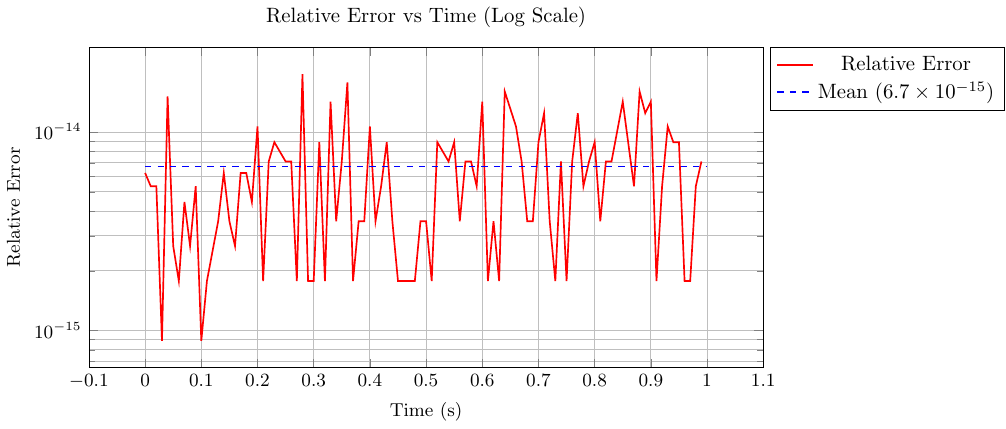}
    \caption{Validation of derivative transfer via spline representation, relative $L^2$-error of the heat flux remains at machine precision level.}
    \label{fig:relative-error-flux}
\end{figure}
\FloatBarrier

As illustrated in Figure~\ref{fig:relative-error-flux}, the spline-based temperature field received from the $\Omega_D$ participant on the interface $\Gamma_D$ is differentiated to evaluate the heat flux. These computed flux values are then compared against the exact analytical solution. The red curve presents the relative error at each time step on a logarithmic scale, while the blue dashed line denotes the mean error across all time steps, which remains consistently around $6.7 \times 10^{-15}$. This result demonstrates the high-fidelity derivative that can be accurately evaluated from the spline field without needing to transfer the derivative separately.

\subsection{Perpendicular flap}
\begin{figure}[h!]
    \centering
    \begin{minipage}{0.60\textwidth}
        \centering
        \includegraphics[width=\textwidth]{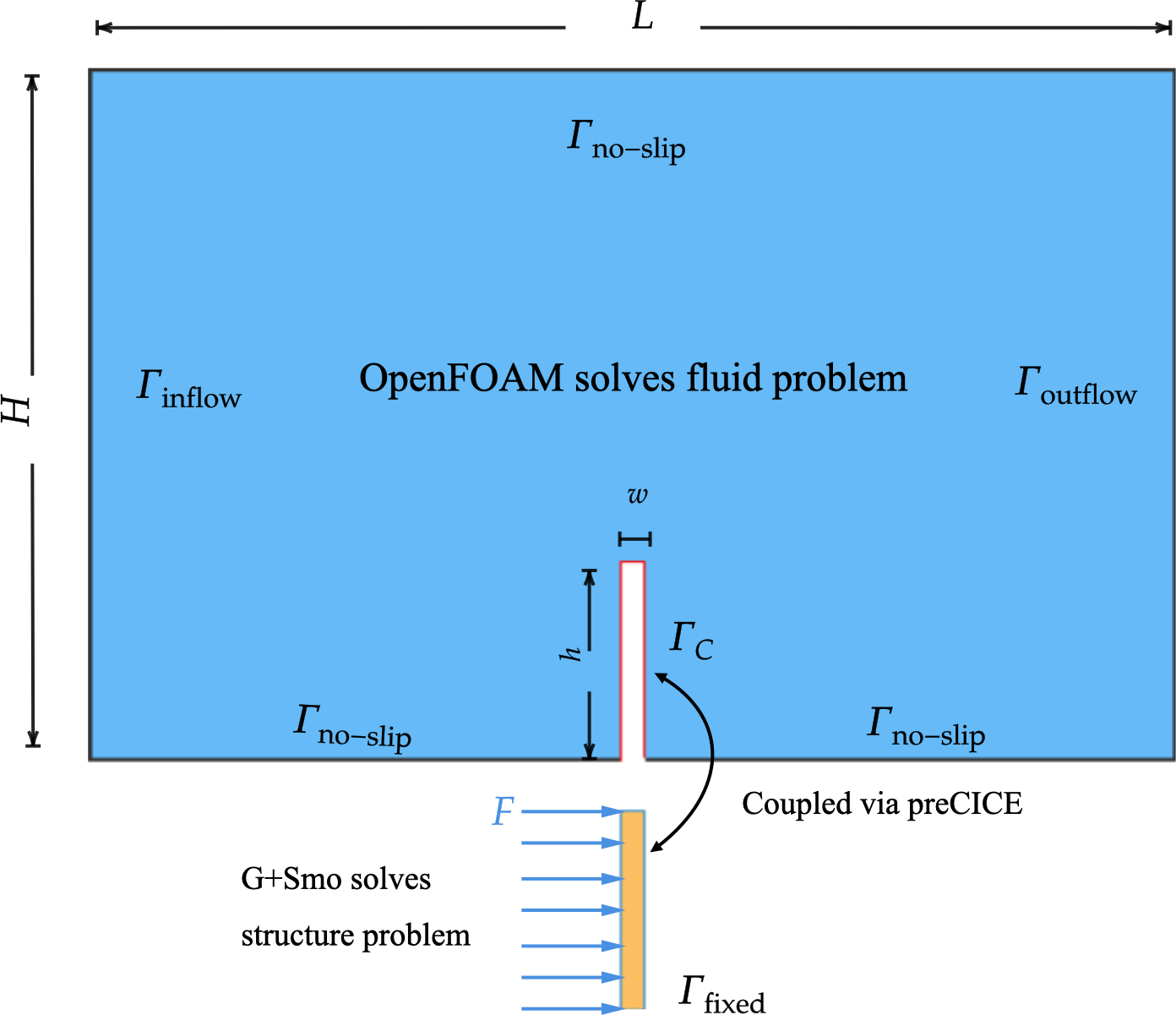}  
    \end{minipage}%
    \hfill
    \begin{minipage}{0.3\textwidth}
        \centering
        \renewcommand{\arraystretch}{1.2}  
        \begin{tabular}{@{}l l l@{}}
            \toprule
            \multicolumn{3}{c}{\textbf{Dimensions}} \\
            \midrule
            $L$ & 6.00 & m \\
            $H$ & 4.00 & m \\
            $l$ & 2.95 & m \\
            $h$ & 1.00 & m \\
            $w$ & 0.10 & m \\
            \midrule
            \multicolumn{3}{c}{\textbf{Fluid}} \\
            \midrule
            $\mathrm{Ma}_\infty$ & 0.01 & \\
            $p_\infty$ & 101325 & Pa \\
            $T_\infty$ & 288.15 & K \\
            \midrule
            \multicolumn{3}{c}{\textbf{Solid}} \\
            \midrule
            $E$ & $4 \cdot 10^6$ & N/m$^2$ \\
            $\nu_s$ & $3 \cdot 10^{-1}$ & \\
            $\rho_s$ & $3 \cdot 10^3$ & kg/m$^3$ \\
            \bottomrule
        \end{tabular}
    \end{minipage}
    \caption{Setup and parameters of the perpendicular flap test case. The original image is taken from the preCICE paper \cite{Chourdakis2021}.}
    \label{figure:perpendicular-flap-setup}
\end{figure}

The perpendicular flap benchmark is a widely used FSI test case \cite{Chourdakis2021} provided by preCICE to evaluate the performance and accuracy of partitioned coupling strategies between fluid and solid solvers. This benchmark is particularly useful for assessing solver interoperability and conducting a mesh convergence study. 
\FloatBarrier
We model a three-dimensional flow that interacts with an elastic flap fixed perpendicularly to the channel floor. The physical parameters for both the fluid and solid domains are detailed in Figure~\ref{figure:perpendicular-flap-setup}.
\FloatBarrier
Figure~\ref{fig:perpendicular-flap-displacement} illustrates the time evolution of the flap tip displacement under fluid loading with time step size.  The flap exhibits characteristic oscillatory behavior due to the interplay between fluid pressure and structural elasticity. 
\FloatBarrier

\FloatBarrier
\begin{figure}[h!]
    \centering
    \begin{minipage}{0.8\linewidth}
    \centering
    \begin{subfigure}[t]{\linewidth}
        \centering
        \includegraphics[width=0.6\linewidth,trim={300 0 310 1050},clip]{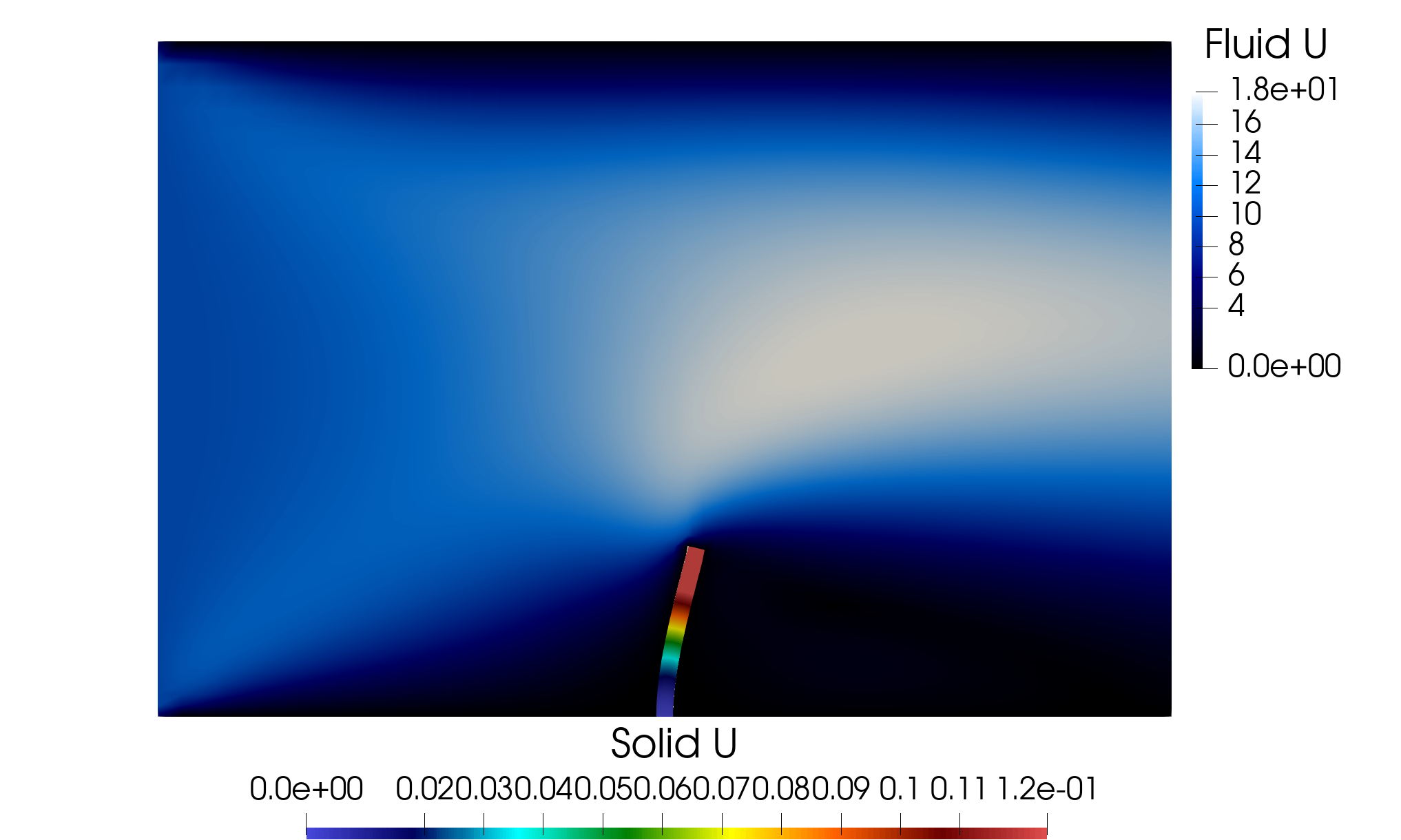}
    \end{subfigure}
    
    \begin{subfigure}[t]{0.49\linewidth}
        \centering
        \includegraphics[width=\linewidth,trim={300 170 310 0},clip]{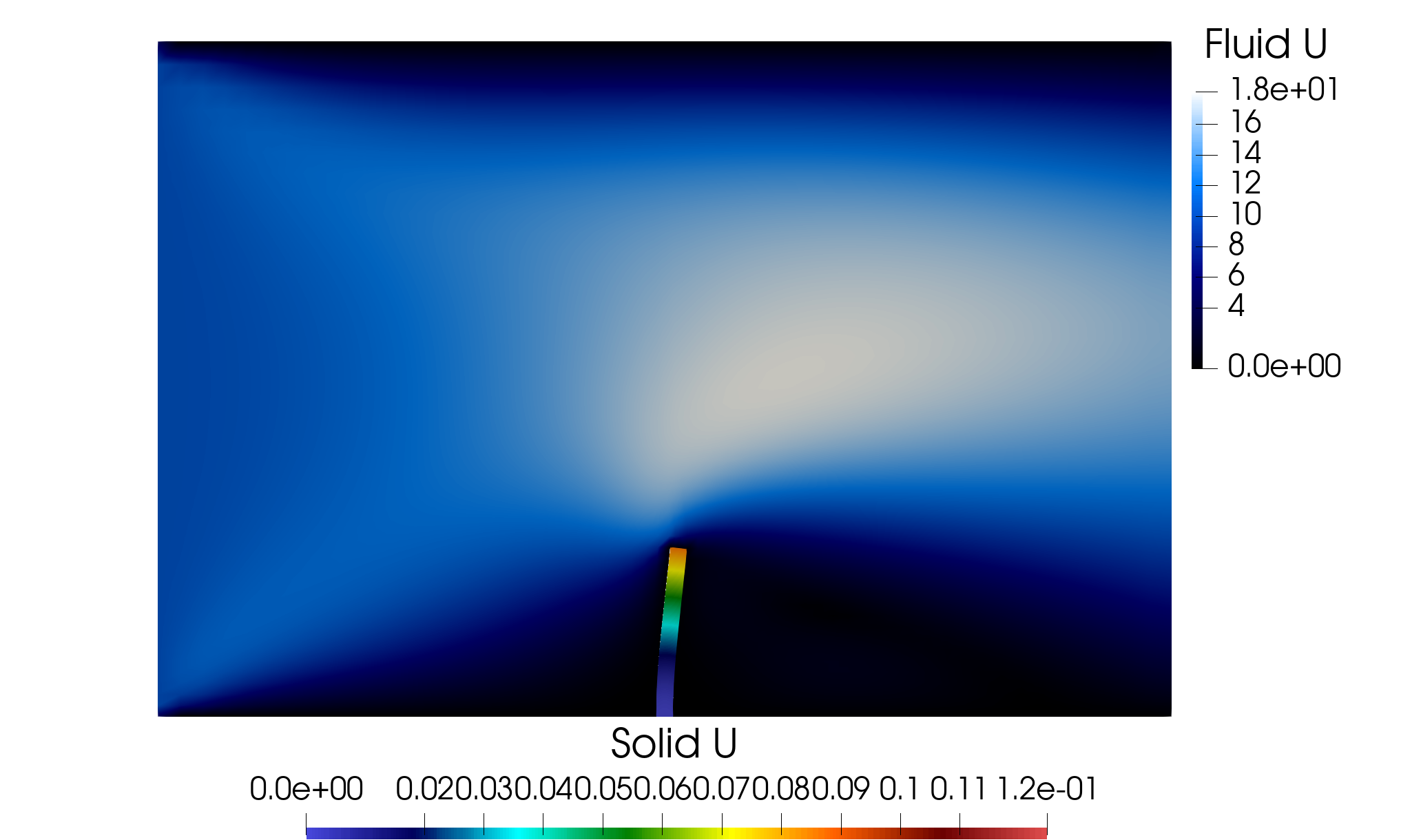}
        \caption{$t=0.3\,\text{s}$}
        \label{fig:perpendicular-flap-displacement0.3}
    \end{subfigure}
    \hfill
    \begin{subfigure}[t]{0.49\linewidth}
        \centering
        \includegraphics[width=\linewidth,trim={300 170 310 0},clip]{figures/perpendicular-flap0.7.png}
        \caption{$t=0.7\,\text{s}$}
        \label{fig:perpendicular-flap-displacement0.7}
    \end{subfigure}
    \end{minipage}
    \hfill
    \begin{minipage}{0.12\linewidth}
        \centering
        \includegraphics[width=\linewidth,trim={1710 600 0 0},clip]{figures/perpendicular-flap0.7.png}
    \end{minipage}
    \caption{Coupled fluid–structure interaction at selected time instances: fluid velocity magnitude (background) and structural pressure distribution (foreground) at $t = 0.3\,\text{s}$ and $t = 0.7\,\text{s}$.}
    \label{fig:perpendicular-flap-displacement}
\end{figure}
In Figure~\ref{fig:perpendicular-flap-displacement}, the dark blue curve corresponds to G+Smo, which is the only solver in this comparison based on IGA. It shows strong agreement with other well-established solvers in both the amplitude and phase of the oscillatory motion of the flap tip. This comparison is based on the tip displacement in the x-direction, using a time step size of $\Delta t = 0.01\,\text{s}$.
\FloatBarrier
\begin{figure}[htbp]
    \centering

    \begin{minipage}[t]{0.45\linewidth}
        \centering
        \begin{tikzpicture}
          \begin{axis}[
            width=\linewidth, 
            height=8cm,
            xlabel={Time (s)},
            ylabel={Tip displacement in \(x\) (m)},
            legend style={
              legend columns=3,
              at={(0.5,-0.2)},
              anchor=north,
            },
            cycle list={
              {gray!60!green, thick, solid},            
              {gray!30!red, thick, dashed},             
              {gray!50!teal, thick, dotted},            
              {gray!30!orange, thick, dash dot},        
              {gray!30!purple, thick, densely dotted},  
              {blue!70!gray, very thick, solid},        
            },
          ]
        
            \addplot+[] table[col sep=space, header=true, x=Time, y=Displacement0]
              {doc/openfoam-calculix-v2404.log};
            \addlegendentry{CalculiX}
        
            \addplot+[] table[col sep=space, header=true, x=Time, y=Displacement0]
              {doc/openfoam-dealii-v2404.log};
            \addlegendentry{deal.II}
        
            \addplot+[] table[col sep=space, header=true, x=Time, y=Displacement0]
              {doc/openfoam-dune-v2404.log};
            \addlegendentry{DUNE}
        
            \addplot+[] table[col sep=space, header=true, x=Time, y=Displacement0]
              {doc/openfoam-fenics-v2404.log};
            \addlegendentry{FEniCS}
        
            \addplot+[] table[col sep=space, header=true, x=Time, y=Displacement0]
              {doc/openfoam-nutils-v2404.log};
            \addlegendentry{Nutils}
        
            \addplot+[] table[col sep=space, header=true, x=Time, y=Displacement0]
              {doc/openfoam-gismo-v2404.log};
            \addlegendentry{\textbf{G+Smo}}
        
          \end{axis}
        \end{tikzpicture}
        \caption{Displacement of the perpendicular flap over time under fluid loading.}
        \label{fig:perpendicular-flap-displacement}
    \end{minipage}
    \hfill
    \begin{minipage}[t]{0.45\linewidth}
        \centering
        \begin{tikzpicture}
        \begin{loglogaxis}[
            width=\linewidth,
            height=8cm,
            xlabel={Element size $h$ ($=2^{-r}$)},
            ylabel={$L^2$ Error in $x$-displacement},
            log basis x=2,
            log basis y=10,
            xtick={0.5, 0.25, 0.125},
            xticklabels={$2^{-1}$, $2^{-2}$, $2^{-3}$},
            legend style={
              at={(0.5,-0.2)},
              anchor=north,
            },
        ]
        
        \addplot+[mark=o, solid, thick] coordinates {
            (0.5, 7.337529e-03)
            (0.25, 1.984676e-03)
            (0.125, 5.399622e-04)
        };
        \addlegendentry{Element order $p=3$}
        
        \addplot+[mark=square*, solid, thick] coordinates {
            (0.5, 6.774220e-04)
            (0.25, 3.399797e-04)
            (0.125, 1.188624e-04)
        };
        \addlegendentry{Element order $p=4$}
        
        \end{loglogaxis}
        \end{tikzpicture}
        \caption{Convergence analysis of flap tip displacement.}
        \label{fig:convergence-analysis}
    \end{minipage}
\end{figure}
\FloatBarrier
A mesh convergence study has been conducted in Figure~\ref{fig:convergence-analysis}. The $L^2$ error in the x-displacement is evaluated over the entire simulation time for different mesh sizes and element orders. Two cases are shown: element order $p = 3$ (blue dashed line) and $p = 4$ (red solid line with square markers).

The error is computed by comparing each solution to a reference solution obtained using a fine mesh with $r = 4$ uniform refinements and a high-order spline basis with $p = 5$. The horizontal axis represents the element size $h = 2^{-r}$, where $r$ is the refinement level. The results demonstrate that increasing the element order and refinement level yields solutions that converge toward the reference solution.


\subsection{Turek Hron FSI Benchmark}
The third benchmark involves two-dimensional flow-induced oscillations of an elastic structure mounted behind a rigid cylinder, commonly referred to as the Turek-Hron FSI3 benchmark. This configuration is widely used to evaluate the performance of FSI solvers, particularly under challenging conditions characterized by strong coupling effects and significant added-mass phenomena. 

\begin{figure}[h!]
    \centering
    \includegraphics[width=0.8\linewidth]{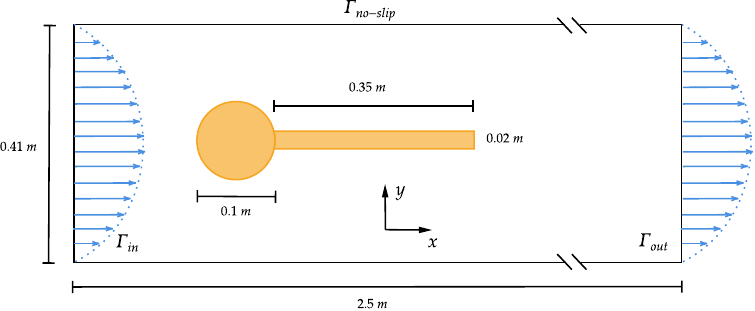}
    \caption{Turek Hron FSI benchmark setup.}
    \label{fig:turek-hron-setup}
\end{figure}

The geometry of the Turek-Hron FSI3 benchmark is illustrated in Figure~\ref{fig:turek-hron-setup}. The configuration consists of a flexible beam attached to the downstream side of a fixed circular cylinder, located inside a rectangular channel of height $0.41\,\mathrm{m}$ and length $2.5\,\mathrm{m}$. The diameter of the cylinder is $0.1\,\mathrm{m}$, and the beam has a length of $0.35\,\mathrm{m}$ and a thickness of $0.02\,\mathrm{m}$. The inflow boundary $\Gamma_{\text{in}}$ is prescribed with a parabolic velocity profile, $v_{\text{in}}(y) = 1.5 \bar{v} \frac{4y(0.41 - y)}{0.41^2}$, where $\bar{v}$ is the average inlet velocity. No-slip boundary conditions are applied along the top and bottom channel walls, denoted by $\Gamma_{\text{no-slip}}$. The average inlet velocity $\bar{v} = 2\,m/s$.

\FloatBarrier
\FloatBarrier

The fluid is assumed to be Newtonian and incompressible, with a density of $\rho_f = 1000\,\mathrm{kg/m^3}$ and a kinematic viscosity of $\nu_f = 10^{-4}\,\mathrm{m^2/s}$, resulting in a Reynolds number of $\mathrm{Re} = 200$. The fluid domain is discretized using an OpenFOAM structured mesh composed mainly of hexahedral cells, consisting of approximately $46$k cells.  The structure is modelled using the Saint-Venant–Kirchhoff material with a Young's modulus of $E = 5.6\,\mathrm{MPa}$ and Poisson's ratio $\nu_s = 0.4$. The solid density is taken as $\rho_s = 1000\,\mathrm{kg/m^3}$, which matches the fluid density to ensure neutral buoyancy.

\FloatBarrier
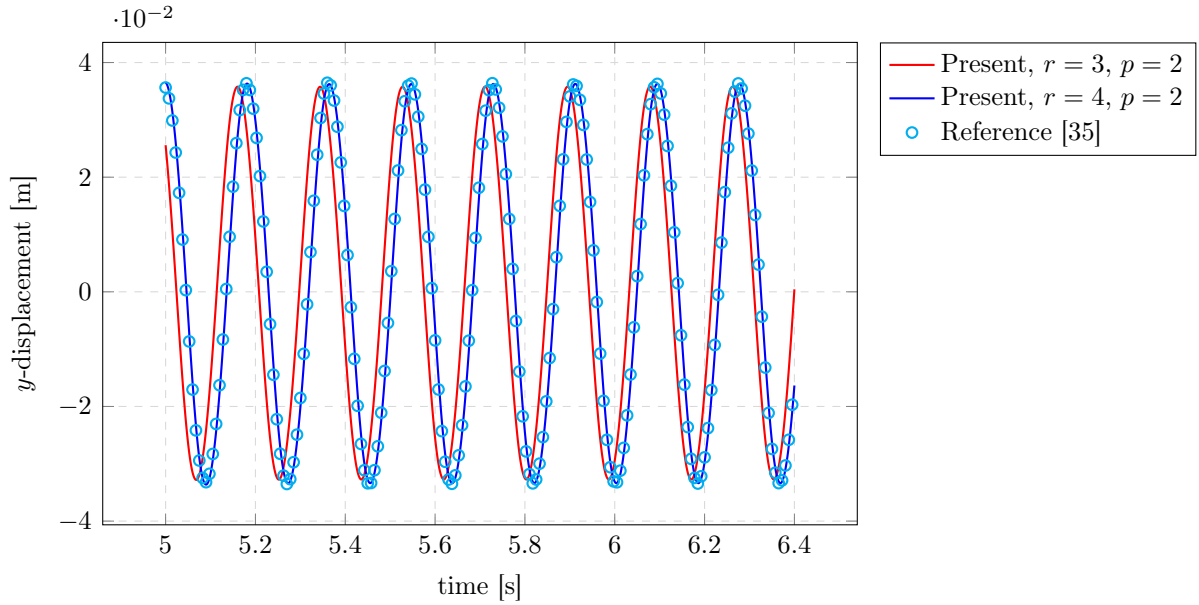
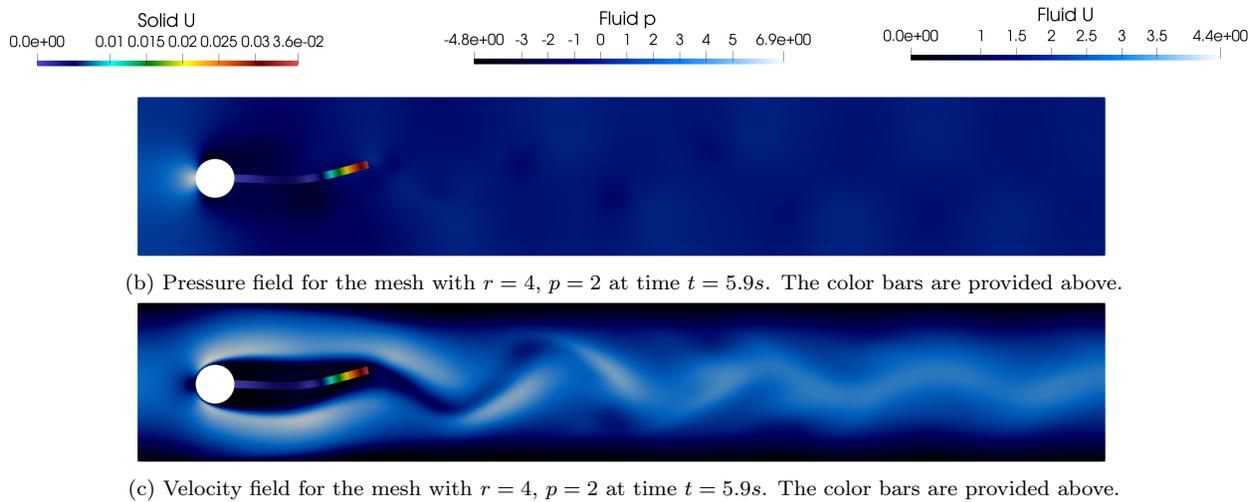
\begin{figure}[h!]
\centering

\begin{subfigure}[t]{\linewidth}
\begin{tikzpicture}
\begin{axis}
[
    width=0.7\linewidth,
    height=8cm,
    xlabel={time [s]},
    ylabel={$y$-displacement [m]},
    restrict x to domain = 5.0:6.4,
    legend pos = outer north east,
    cycle list name=customcolorlist
]

\addplot+[no markers] table [x index=0, y index=4, col sep=space] {validation_data/precice-Solid-watchpoint-Flap-Tip-r3-run.log};
\addlegendentry{Present, $r=3$, $p=2$}
\pgfplotsset{cycle list shift=1}

\addplot+[no markers] table [x index=0, y index=4, col sep=space] {validation_data/precice-Solid-watchpoint-Flap-Tip-r4-run.log};
\addlegendentry{Present, $r=4$, $p=2$}
\pgfplotsset{cycle list shift=2}

\addplot+[mark=o,cyan,only marks,each nth point=30] table [x index=0, y index=11, col sep=space] {validation_data/ref_fsi3.point};
\addlegendentry{Reference \cite{featflow_fsi3}}
\end{axis}
\end{tikzpicture}
\caption{$y$-displacement of the flap tip with respect to time on the interval $t\in[5.0,6.4]$. The reference data \cite{featflow_fsi3} is represented by blue marker, and the data obtained by the present model is represented by blue and red lines.  }
\label{fig:turek-hron-fsi3-a}
\end{subfigure}

\begin{subfigure}[t]{0.3\linewidth}
\centering
\includegraphics[width=\linewidth,trim={1500 105 100 1010},clip]{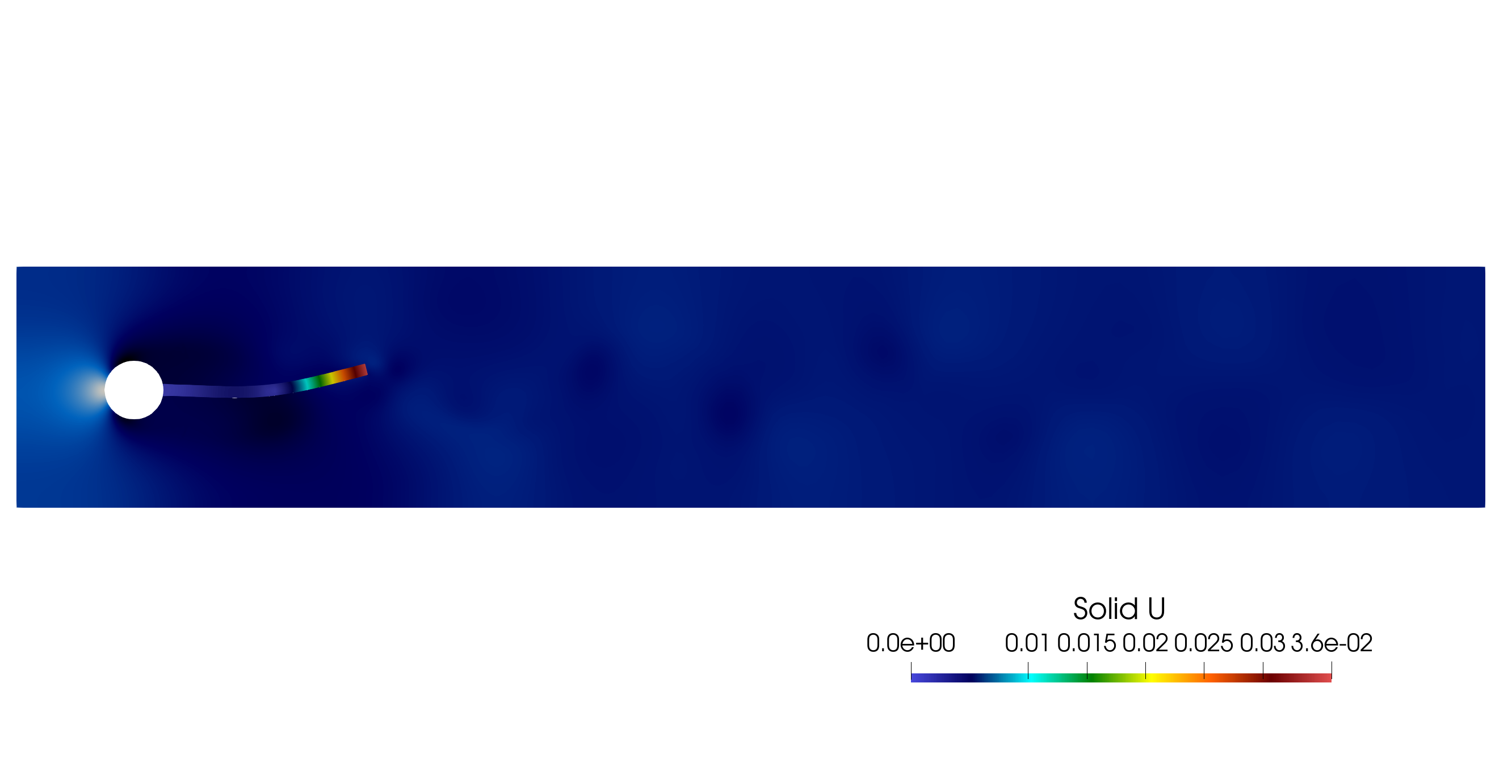}
\end{subfigure}
\hfill
\begin{subfigure}[t]{0.3\linewidth}
\centering
\includegraphics[width=\linewidth,trim={1500 100 100 1010},clip]{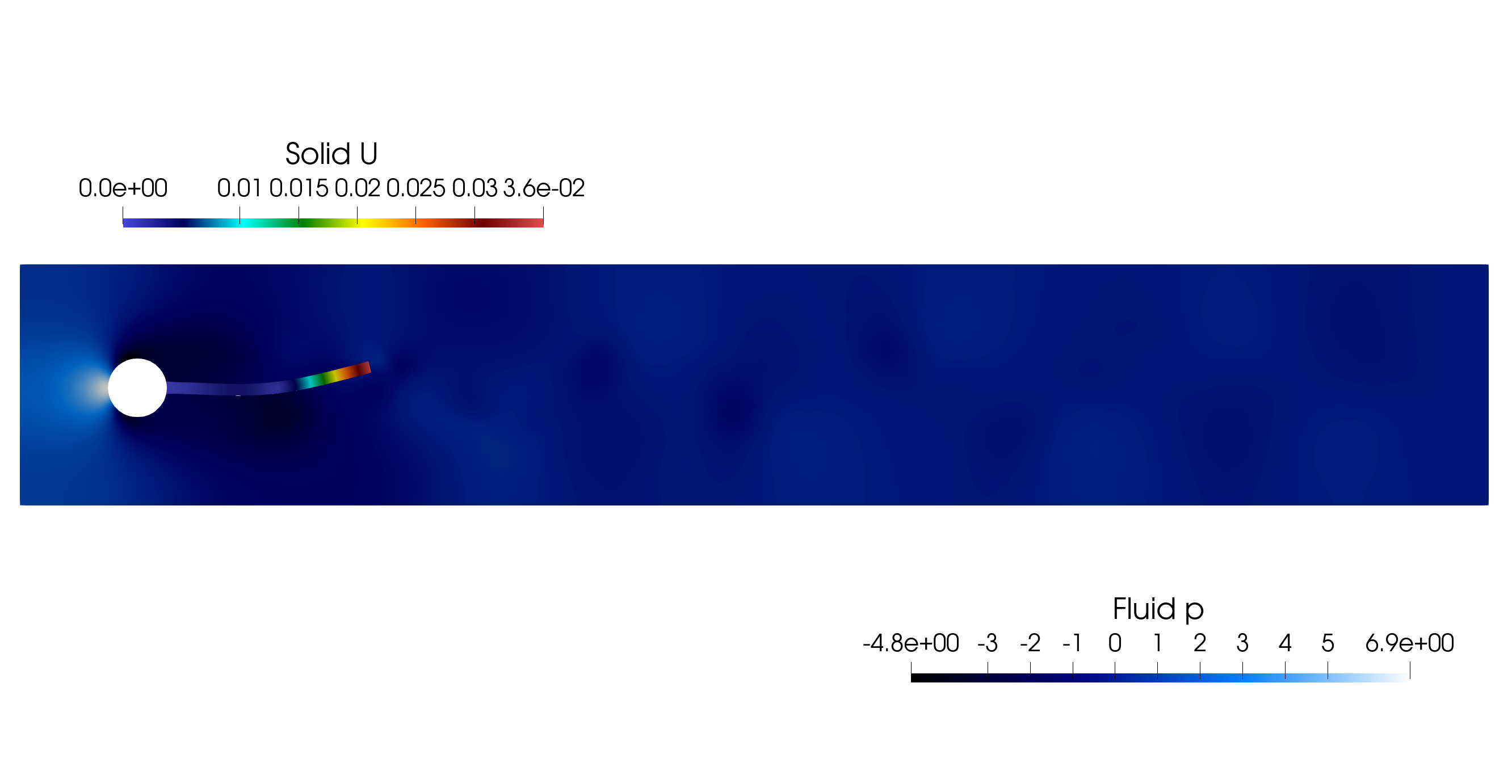}
\end{subfigure}
\hfill
\begin{subfigure}[t]{0.3\linewidth}
\centering
\includegraphics[width=\linewidth,trim={1500 90 100 1010},clip]{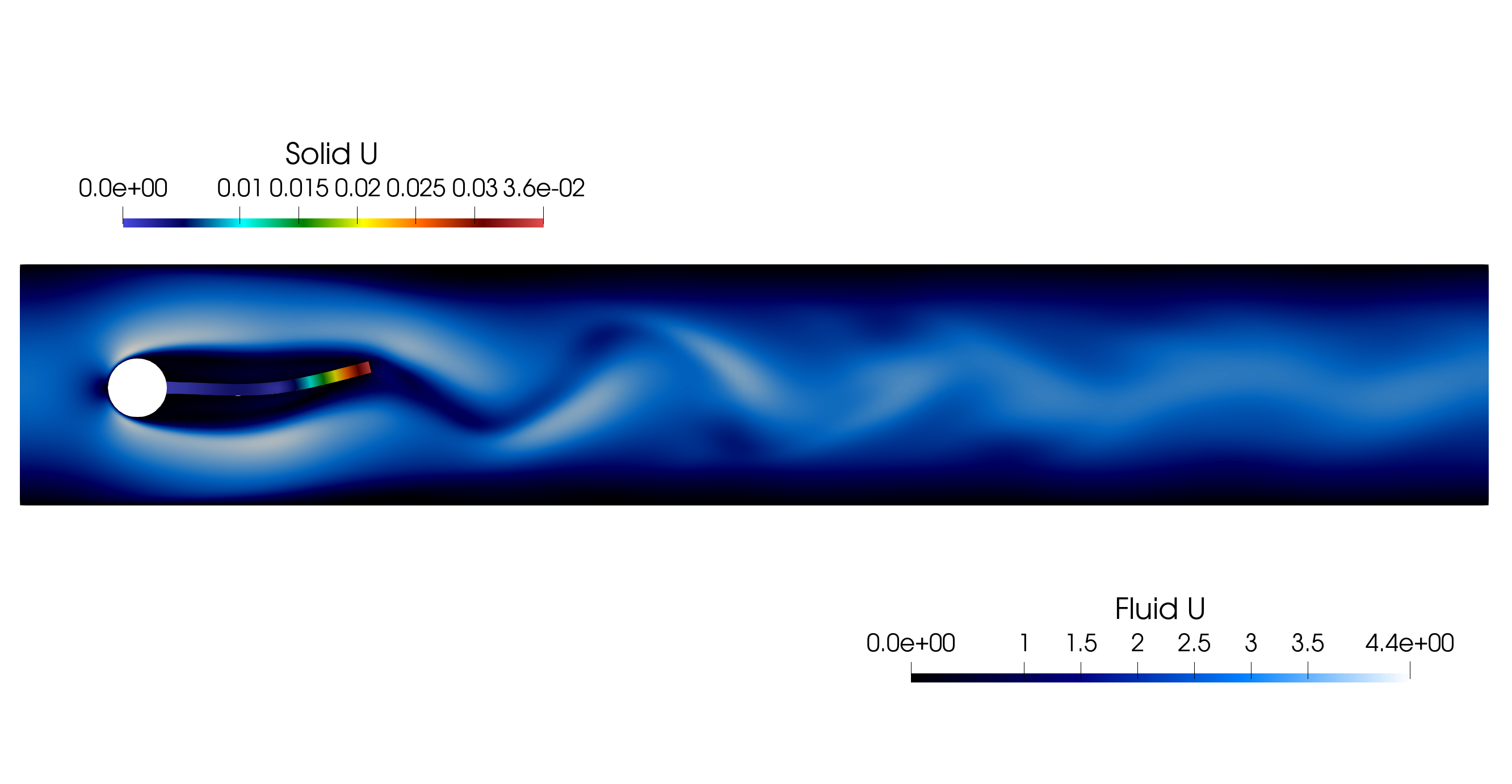}
\end{subfigure}

\begin{subfigure}[t]{\linewidth}
    \centering
    \includegraphics[width=\linewidth,trim={0 430 0 430},clip]{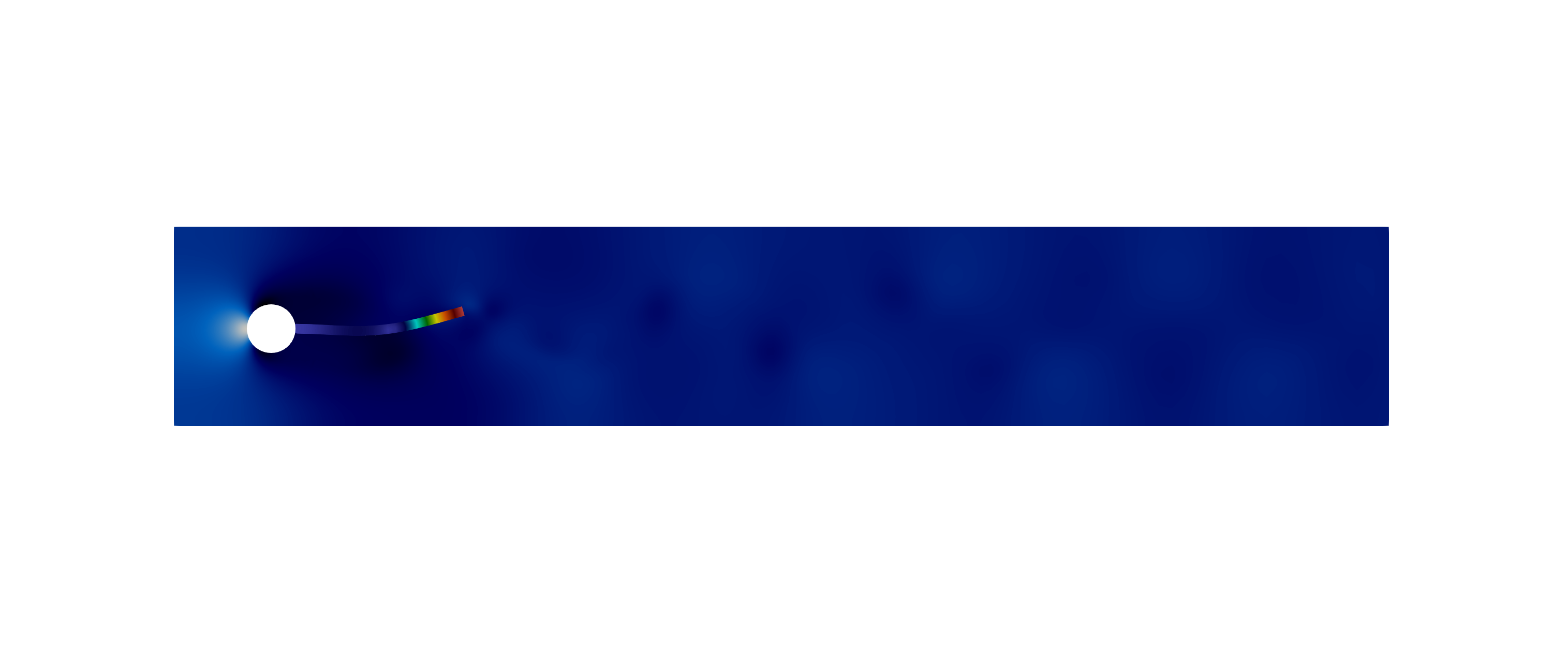}
    \caption{Pressure field for the mesh with $r=4$, $p=2$ at time $t=5.9s$. The color bars are provided above.}
    \label{fig:turek-hron-fsi3-b}
\end{subfigure}

\begin{subfigure}[t]{\linewidth}
    \centering
    \includegraphics[width=\linewidth,trim={0 430 0 430},clip]{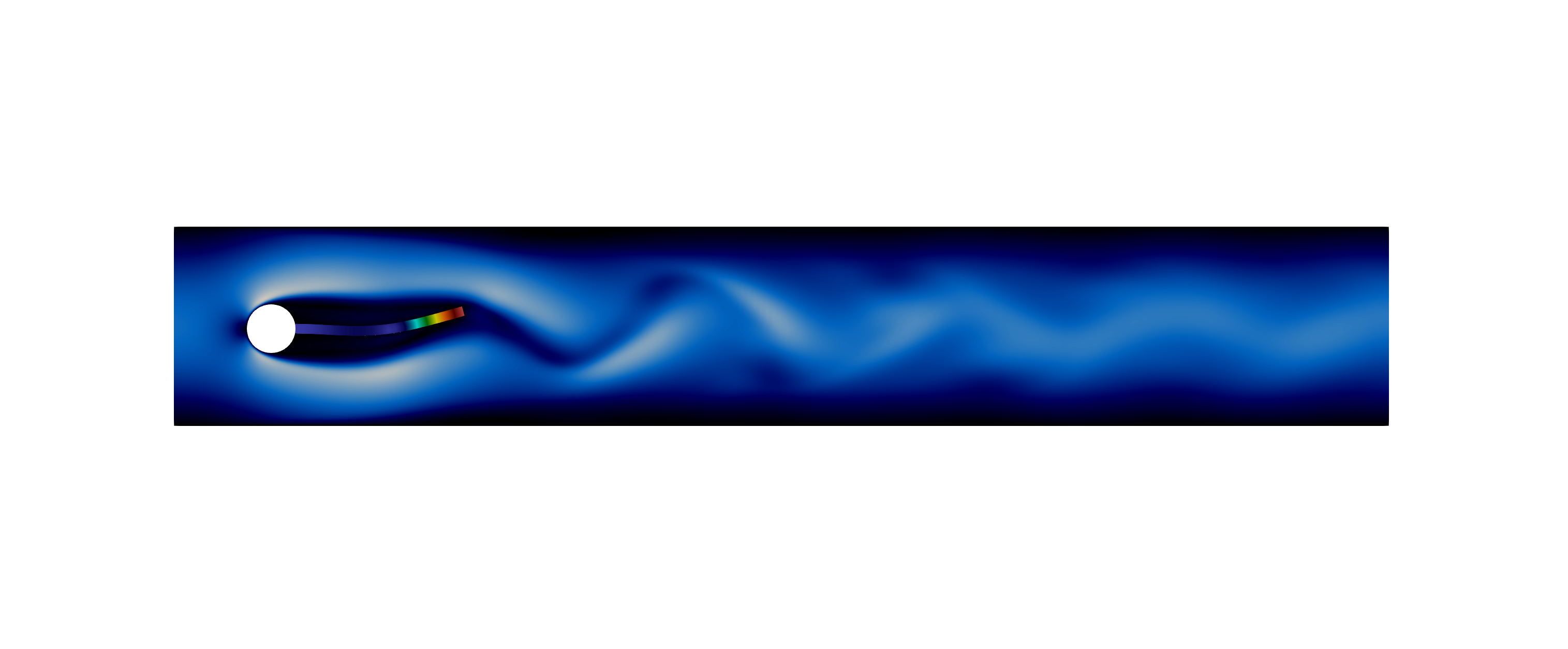}
    \caption{Velocity field for the mesh with $r=4$, $p=2$ at time $t=5.9s$. The color bars are provided above.}
    \label{fig:turek-hron-fsi3-c}
\end{subfigure}
\caption{Simulation results for the Turek-Hron FSI3 benchmark using G+Smo coupled with OpenFOAM. (a) Vertical displacement of the beam tip over time. (b) Beam displacement and fluid pressure distribution. (c) Beam displacement and fluid velocity magnitude for mesh $r=4$ with degree $p=2$ at time step $t=5.9 \text{ s}$. 
}
\label{fig:turek-hron-fsi3-result}
\end{figure}

The left plot in Figure~\ref{fig:turek-hron-fsi3-result} shows the time evolution of the vertical displacement at the tip of the elastic beam. The results obtained using G+Smo with two refinement levels ($r = 3$ and $r = 4$, both with degree $p = 2$) are compared against reference data. The numerical solution from the finer mesh ($r = 4$) closely follows the reference in both amplitude and phase, while the coarser mesh ($r = 3$) slightly underestimates the peak displacements but still captures the oscillation frequency well.

This strong agreement demonstrates that our coupling framework for the IGA-based structure solver is capable of accurately resolving complex FSI dynamics. It highlights the method's ability to perform high-fidelity and efficient coupling within libraries that support splines, while also remaining compatible with solvers that work without spline-based boundary representations.

\FloatBarrier
\section{Conclusion}
\label{section:conclusion}
In this work, we have put together a seamless, start-to-finish spline coupling workflow for partitioned multiphysics simulations, with a particular focus on fluid-structure interaction problems. 
Our approach demonstrates that splines not only serve as mathematical representations for geometry and solution fields, but also function effectively as an efficient data exchange format between different solvers.

\begin{enumerate}
    \item Development of spline-based coupling methods for enhanced representation of interface fields. 
    \item Comprehensive analysis of communication overhead, demonstrating the efficiency of spline-based methods over traditional vertex-based approaches. Our theoretical model and experimental validation show that communication overhead for vertex-based coupling grows rapidly with mesh refinement, while spline-based methods demonstrate significantly lower growth rates. 
    \item Implementation of robust coupling strategies that naturally handle non-conforming but nested meshes and preserve solution derivatives across interfaces. The spline-based approach eliminates the need for complex interpolation schemes when meshes do not match at the interface, and inherently preserves both solution values and their gradients, as demonstrated in our heat conduction benchmark where derivative information was accurately maintained without requiring separate transfer operations.
\end{enumerate}

The proposed methods are verified using vertical beam test cases under both constant and nonlinear loading, and validated through established FSI benchmarks, including the perpendicular flap and Turek-Hron configurations. Although all benchmarks are two-dimensional, this restriction does not limit the generality of our findings. The primary goal is to assess the coupling strategy and its communication behavior—both of which extend naturally to three-dimensional settings. 

Future work will focus on extending the proposed methods to more complex multiphysics problems, including those involving large deformations, topology changes, and multiple coupled fields. 
We plan to implement adaptive local refinement using Truncated Hierarchical B-splines (THB-splines) \cite{GIANNELLI2016337}, which enable targeted mesh refinement in regions with high solution gradients while preserving the advantages of spline-based coupling in 3D geometries. Additionally, we intend to adapt our coupling methodology to Kirchhoff--Love shell structures \cite{Hugo_thesis,KIENDL20093902}, which are prevalent in many FSI applications such as aerospace and biomedical engineering. Another promising direction is the development of IGA-based mesh transfer techniques specifically optimized for FSI applications, which could further improve the handling of moving and deforming interfaces \cite{JI2021113615}. 
Finally, we will continue optimizing our implementation to support large-scale, industrial-grade simulations.

\section*{CRediT authorship contribution statement}
\textbf{J. Li:} Writing – original draft, Visualization, Validation, Software, Methodology, Investigation, Formal analysis, Conceptualization. \textbf{H.M. Verhelst:} Writing – review \& editing, Software, Investigation, Conceptualization, Supervision. \textbf{J.H. Den Besten:} Writing – review, Supervision, Project administration, Funding acquisition. \textbf{M. M\"{o}ller:} Writing – review, Software, Supervision, Project administration, Funding acquisition, Conceptualization.

\section*{Declaration of competing interest}
The authors declare that they have no known competing financial interests or personal relationships that could have appeared to influence the work reported in this paper.

\section*{Acknowledgments}
We would like to express our gratitude to Jun.-Prof. Benjamin Uekermann and Gerasimos Chourdakis from the preCICE team for their technical support and valuable insights regarding the coupling interface implementation. 
We also appreciate Dr. Marin Lauber for the fruitful discussions on integrating spline-based coupling methods in Waterlily.jl. 

This research has been conducted as a part of the project FlexFloat with project number 19002 of the research programme Open Technology which is (partly) financed by the Dutch Research Council (NWO), Netherlands. In addition, HMV is grateful for the financial support from the Italian Ministry of University and Research (MUR) through the PRIN projects COSMIC (no. 2022A79M75) and ASTICE (no. 202292JW3F), with the contribution of the European union – Next Generation EU.

\section*{Data availability}
The proposed coupling methods are implemented as part of the \href{https://github.com/gismo/gismo}{\texttt{\textcolor{blue}{G+Smo}}} (Geometry + Simulation Modules) repository. We used \href{https://precice.org}{\texttt{\textcolor{blue}{preCICE}}} to interface with the fluid solver OpenFOAM and the structure solver G+Smo.
\bibliographystyle{plainnat} 
\bibliography{reference}

\end{document}